\documentclass[dvips,ejs,preprint]{imsart}

\usepackage{graphicx,bbm,amssymb,amsmath,amsfonts,amsthm,stmaryrd}
\usepackage[pdftex]{hyperref}
\usepackage[numbers]{natbib}

\RequirePackage[numbers]{natbib}
\RequirePackage[colorlinks,citecolor=blue,urlcolor=blue]{hyperref}

\RequirePackage{hypernat}

\def\R{\mathbb{R}}

\def\N{\mathbb{N}}
\def\Z{\mathbb{Z}}
\def\E{\mathbb{E}}

\def\Var{\textrm{Var}}

\def\F{\mathcal{F}}

\def\x{\mathcal{X}}
\def\G{\mathcal{G}}
\def\Ph{\widehat{P}}

\def\Vh{\widehat{V}}
\def\Nh{\widehat{N}}
\def\ph{\widehat{p}}
\def\mh{\widehat{m}}

\def\Gh{\hat{G}}

\DeclareMathOperator*{\argmin}{arg\,min}
\DeclareMathOperator*{\argmax}{arg\,max}
\def\pen{\mathrm{pen}}

\newcommand{\paren}[1]{\left( \left. #1 \right. \right)} 
 
\newcommand{\set}[1]{\left\{ \left. #1 \right. \right\}}
\newcommand{\absj}[1]{\left\lvert #1 \right\rvert} 
\providecommand{\norm}[1]{\left \lVert #1 \right\rVert}

\theoremstyle{plain}
\newtheorem{theorem}{Theorem}[section]
\newtheorem{lemma}[theorem]{Lemma}
\newtheorem{proposition}[theorem]{Proposition}
\newtheorem{coro}[theorem]{Corollary}
\newtheorem{definition}[theorem]{Definition}

\begin{document}

\begin{frontmatter}

\title{An Oracle Approach \\
for Interaction Neighborhood Estimation in Random Fields}
\runtitle{Interaction Neighborhood Estimation}


\author{\fnms{Matthieu} \snm{Lerasle}\thanksref{t1}\ead[label=e1]{lerasle@gmail.com}}
\and
\author{\fnms{Daniel, Y.} \snm{Takahashi}\thanksref{t2}\ead[label=e2]{takahashiyd@gmail.com}}
\thankstext{t1}{Supported by FAPESP grant 2009/09494-0.} 
\address{Instituto de Matem\'atica e Estat\'{\i}stica\\
Universidade de S\~ao Paulo\\
Caixa Postal 66281\\
05315-970 S\~ao Paulo, Brasil\\
\printead{e1}}

\thankstext{t2}{Supported by FAPESP grant 2008/08171-0. }
\address{Instituto de Matem\'atica e Estat\'{\i}stica\\
Universidade de S\~ao Paulo\\
Caixa Postal 66281\\
05315-970 S\~ao Paulo, Brasil\\
currently at Princeton University\\
Department of Psychology and Neuroscience\\
Green Hall, Princeton, NJ 08648\\
\printead{e2}}

\runauthor{M. Lerasle and D. Y. Takahashi}

\begin{abstract}
We consider the problem of interaction neighborhood estimation from the partial observation of a finite number of realizations of a random field. We introduce a model selection rule to choose estimators of conditional probabilities among natural candidates. Our main result is an oracle inequality satisfied by the resulting estimator. We use then this selection rule in a two-step procedure to evaluate the interacting neighborhoods. The selection rule selects a small prior set of possible interacting points and a cutting step remove from this prior set the irrelevant points.\\
We also prove that the Ising models satisfy the assumptions of the main theorems, without restrictions on the temperature, on the structure of the interacting graph or on the range of the interactions. It provides therefore a large class of applications for our results. We give a computationally efficient procedure in these models. We finally show the practical efficiency of our approach in a simulation study.
\end{abstract}





\end{frontmatter}

\section{Introduction}

\indent 
Graphical models, also known as random fields, are used in a variety of domains, including computer vision \citep{Besag93, Woods78}, image processing \citep{Cross83}, neuroscience \citep{Schneidman06}, and as a general model in spatial statistics \citep{Ripley81}. The main motivation for our work comes from neuroscience where the advancement of multichannel and optical technology enabled the scientists to study not only a unit of neurons per time, but tens to thousands of neurons simultaneously \citep{Takahashi10}. The very important question now in neuroscience is to understand how the neurons in this ensemble interact with each other and how this is related to the animal behavior \citep{Schneidman06, Brown04}. This question turns out to be hard for three reasons at least. First, the experimenter has always only access to a small part of the neural system. Moreover, there is no really good model for population of neurons in spite of the good models available for single neurons. Finally, strong long range interactions exist \citep{Li10}. Our work tries to overcome some of these difficulties as will be shown.

A random field can be specified by a discrete set of sites $G$, possibly infinite, a finite alphabet of spins $A$, and a probability measure $P$ on the set of configurations $\x(G)=A^G$. One of the objects of interest are the one-point specification probabilities, defined for all sites $i$ in $G$ and all configurations $x$ in $\x(G)$ by a regular version of the conditional probability
$$P\left(\; x(i)\;  | \; x(j),\; j\in G/\{i\}\right).$$ 
From a statistical point of view, two problems are of natural interest. 

\vspace{0.2cm}

\noindent
{\bf Interaction neighborhood identification problem (INI):} 

The INI problem is to identify, for all sites $i$ in $G$, the minimal subset $G_i$ of $G$ necessary to describe the specification probabilities in site $i$ (see Sections \ref{Section:Preliminaries} and \ref{Section:GeneralResults} for details). $G_i$ is called the interaction neighborhood of $i$ and the points in $G_i$ are said to interact with $i$. $G_i$ is not necessarily finite but only a finite subset $V_M\subset G$ of sites is observed. The observation set is a sample $X_{1:n}(V_M)=(X_1(j),...,X_n(j))_{j\in V_M}$, where $(X_1,...,X_n)$ are i.i.d with common law $P$. The question is then to recover from $X_{1:n}(V_M)$, for all $i$ in $V_M$, the sets $G_i\cap V_M$. 

\vspace{0.2cm}

\noindent{\bf Oracle neighborhood problem (ON):} 

The ON problem is to identify, for all $i$ in $G$, a set $\Gh_i=\Gh_i(X_{1:n}(V_M))$, such that the estimation of the conditional probabilities $P\left(x(i)|x(j),j\in G/\{i\}\right)$ by the empirical conditional probabilities $\Ph(x(i)|x(j),j\in \Gh_i)$ has a minimal risk (see Sections \ref{Section:Preliminaries} and \ref{Section:GeneralResults} for details).  $\Gh_i$ is then said to satisfy an oracle inequality and it is also called oracle. We look for oracles among the subsets of $V_M$ and we consider the $L_{\infty}$-distance between conditional probabilities to measure the risk of the estimators. An oracle is in general smaller than $G_i$ because it should balance approximation properties and parsimony.
 
\vspace{0.2cm}

The literature has mainly been focused in the INI problem, see \cite{Bento09, Bresler08, Csiszar06, GOT10, PW10} for examples. It requires in general strong assumptions on $P$ to be solved. For example, the $\ell_1$-penalization procedure proposed in \cite{PW10} requires an incoherence assumption on the interaction neighborhoods that is very restrictive, as shown by \cite{Bento09}.
Moreover, it is assumed in \cite{Bento09,  Bresler08, PW10} that $G$ is finite and that all the sites are observed, {\it i.e.} $V_M=G$. \citet{Csiszar06} consider the case when $G = \Z^d$ but assume a uniform bound on the cardinality of $G_{i}$. The procedure proposed in \cite{GOT10} holds for infinite graph with each site having infinite neighborhoods, but requires that the main interactions belong to a known neighborhood of $i$ of order $O(\ln n)$. Moreover, the result is proved in the Ising model only when the interaction is sufficiently weak.

The first goal of this paper is to show that the ON problem can be solved without any of these hypotheses. We introduce in Section \ref{subsection:ModelSelection} a model selection criterion to choose a model $\Gh_i$ and prove that it is an oracle in Theorem \ref{theo:GenModSel}. This result does not require any assumption on the structure of the interaction neighborhoods inside or outside $V_M$. 
 
The second objective is to show that a selection rule provides also a useful tool to handle the INI problem. We introduce the following two steps procedure. First, we select, for all sites $i$ in $V_M$, a small subset $\hat{V}_i$ of $V_M$ with the model selection rule. We prove that this set contains the main interacting points inside $V_M$ with large probability. Following the idea introduced in \cite{GOT10}, we use then a test to remove from $\hat{V}_i$ the points of $(G/G_i)\cap \hat{V}_i$. The new test can be applied to all neighborhoods $V_i$ that are smaller than $O(\ln n)$ and that contain the main interaction points in $G_i$. It requires less restrictive assumptions on the interactions outside $V_i$ and on the measure $P$ than the one of \cite{GOT10}. For example, it works in the Ising models without restrictions on the temperature parameter. Furthermore, the two-step method let us look for the interacting points inside all the observation set $V_M$ (of order $O(e^{n^{\beta}})$ for some $0\leq \beta<1$), and not only inside a prior subset $V_i$ (smaller than $O(\ln n)$) of $V_M$. 

All the results hold under a key assumption {\bf H1} that is not classical, but that is satisfied by Ising models, see Theorem \ref{theo:controlBiasGibbs}. We obtain then a large class of models, widely used in practice, where our methods are efficient. We also provide for this model a computationally efficient version of our main algorithms.

The paper is organized as follows. In Section \ref{Section:Preliminaries}, we introduce notations and assumptions used all along the paper. Section \ref{Section:GeneralResults} gives the main results, in a general framework. Section \ref{sect:pairGibbsMeasures} shows the application to Ising models and Section \ref{sect:simu} presents a large simulation study where the problem of the practical calibration of some parameters is adressed. Section \ref{Section:Discussion} is a discussion of the results with a comparison to existing papers. Section \ref{section:Proofs} gives the proofs of the main theorems and some technical results are recalled in an appendix in Section \ref{section:Appendix}. 

\section{Notations and Main Assumptions}\label{Section:Preliminaries}
Let $G$ be a discrete set of \emph{sites}, possibly infinite, $A=\{-1,1\}$ be the binary alphabet of \emph{spins}, and $P$ be a probability measure  on the set of \emph{configurations} $\x(G)=A^G$. More generally, for all subsets $V$ of $G$, let $\x(V)=A^V$ be the set of configurations on $V$. In what follows, the triplet $(G,A,P)$ will be called a \emph{random field}. For all $i$ in $G$, for all $V\subset G$, for all $x$ in $\x(G)$, let $x(V)=(x(j))_{j\in V}$ and for all probability measures $Q$ on $\x(V\cup\{i\})$, let
\begin{equation*}\label{def:probacond}
Q_{i|V}(x)=Q(x(i)|x(V/\{i\}))
\end{equation*}
be a regular version of the conditional probability. All along the paper, we will use the convention that, if $V$ is a finite set, $Q$ a probability measure on $\x(V)$ and $x$ is a configuration such that $Q(x(V/\{i\}))=0$, then $Q_{i|V}(x)=1/2$.\\
For all $x$ in $\x(G)$ and all $j$ in $G$, let $x_j$ be the configuration such that $x_j(k)=x(k)$ for all $k\neq j$ and $x_j(j)=-x(j).$ We say that there is a \emph{pairwise interaction} from $j$ to $i$ if there exists $x$ in $\x(G)$ such that $P_{i|G}(x_{j})\neq P_{i|G}(x).$ For all subsets $V$ of $G$, for all probability measures $Q$ on $\x(V)$, let
$$\omega^{V}_{i,j}(Q)=\sup_{x\in\x(G)}\left\{Q_{i|V}(x)-Q_{i|V}(x_j)\right\}.$$
With the above notations, there is a pairwise interaction from $j$ to $i$ if and only if $\omega^{G}_{i,j}(P)>0$. Our second task in this paper is to recover the set $G_i$ of sites having a pairwise interaction with $i$. This definition differs in general from the one suggested in introduction. However, it is easy to check that they coincide in the Ising models defined in Section \ref{sect:pairGibbsMeasures}.\\
Let $X_{1:n}=(X_1,...,X_n)$ be i.i.d. with common law $P$. Let $V_M$ be a finite subset of $G$ of \emph{observed sites}, with cardinality $M$. The \emph{observation set} is then $X_{1:n}(V_M)=(X_1(V_M),...,X_n(V_M))$. Let $\Ph$ be the \emph{empirical measure} on $\x(G)$ defined for all configurations $x$ in $\x(G)$ by
\begin{equation*}\label{def:est}
\Ph(x)=\frac1n\sum_{i=1}^n\mathbf{1}_{\{X_i(G)=x(G)\}}.
\end{equation*}
For all real valued functions $f$ defined on $\x(G)$, let $\left\|f\right\|_{\infty}=\sup_{x\in \x(G)}|f(x)|$. For all subsets $V$ of $V_M$, the \emph{$L_{\infty}$-risk} of $\Ph_{i|V}$ is defined by $\left\|\Ph_{i|V}-P_{i|G}\right\|_{\infty}$. This risk is naturally decomposed into two terms. From the triangular inequality, we have
$$\left\|\Ph_{i|V}-P_{i|G}\right\|_{\infty}\leq \left\|\Ph_{i|V}-P_{i|V}\right\|_{\infty}+\left\|P_{i|V}-P_{i|G}\right\|_{\infty}.$$
We call variance term the random term $\left\|\Ph_{i|V}-P_{i|V}\right\|_{\infty}$ and bias term the deterministic one $\left\|P_{i|V}-P_{i|G}\right\|_{\infty}$.\\
Let us finally present our general assumptions on the measure $P$. In the following $\nu$ and $\kappa_{\min}$ are positive constants. The two first assumptions are classical and will only be used to discuss the main results.

\vspace{0.2cm}

\noindent
{\bf NN:} {\it (Non-Nullness) For all $x$ in $\x(G)$, $\nu^{-1}\leq P_{i|G}(x)$.}

\vspace{0.2cm}

\noindent
{\bf CA:}{\it (Continuity) For all growing sequences $(V_n)_{n\in\N^*}$ of subsets of $G$ such that  $\cup_{n\in \N^*}V_n=G$, for all $i$ in $G$, }
$$\lim_{n\rightarrow \infty}\left\|P_{i|V_n}-P_{i|G}\right\|_{\infty}=0.$$

\noindent
The following last  assumption is very important for the model selection criterion to work. It is satisfied for example by a generalized form of the Ising model as we will see in Section \ref{sect:pairGibbsMeasures}.
\vspace{0.2cm}

\noindent
{\bf H1:} {\it For all finite subsets $V$ of $G$, }
$$\kappa_{\min}\left\|P_{i|G}-P_{i|V}\right\|_{\infty}\leq \left\|P_{i|G}\right\|_{\infty}-\left\|P_{i|V}\right\|_{\infty}.$$

\section{General results}\label{Section:GeneralResults}
\subsection{Control of the variance term of the $L_{\infty}$-risk}
Our first theorem provides a sharp control of the variance term of the risk of $\Ph_{i|V}$. It holds without assumption on the measure $P$ or the finite subset $V$.
\begin{theorem}\label{theo:controlVar}
Let $P$ be a probability measure on $\x(G)$, let $V$ be a finite subset of $G$. Let $p_-^V=\inf_{x\in \x(G),\;P(x(V))\neq 0}P(x(V))$. There exists an absolut constant $c_1$ such that, for all $\delta>1$,
\begin{equation}\label{eq:ControlVarDet}
P\left(\left\|\Ph_{i|V}-P_{i|V}\right\|_{\infty}> c_1\sqrt{\frac{\ln(\delta/ p^V_-)}{np^V_-}}\right)\leq \frac1{\delta}.
\end{equation}
Moreover, let $\ph_{-}^{V}=n^{-1}\vee\inf_{x\in\x(V)}\Ph(x(V/\{i\}))$. There exists an absolut constant $c_{2}\leq 400$ such that, for all $\delta>1$,
\begin{equation}\label{eq:ControlVarEmp}
P\left(\norm{\Ph_{i|V}(x)-P_{i|V}(x)}_{\infty}> c_2\sqrt{\frac{\ln(\delta n)}{n\ph^V_-}}\right)\leq \frac1{\delta}.
\end{equation}
\end{theorem}

\noindent
{\bf Remark:} 
\begin{itemize}
\item Let $|V|$ denote the cardinality of $V$, if $P$ satisfies {\bf NN} we have $p^V_-\geq \nu^{-|V|}.$ Hence, (\ref{eq:ControlVarDet}) implies that, 
\begin{equation*}
P\paren{\left\|\Ph_{i|V}-P_{i|V}\right\|_{\infty}\leq c_1\sqrt{\nu^{|V|}\frac{|V|\ln(\nu)+2\ln (n)}{n}}}\geq1-n^{-2}.
\end{equation*}
The variance term goes almost surely to $0$ if $\nu^{|V|}< <n(\ln n)^{-1}$. If in addition $P$ satisfies {\bf CA} and $(V_n)_{n\in\N^*}$ is a growing sequence of sets with limit $G$, the estimator $\Ph_{i|V_n}$ is consistent.
\item (\ref{eq:ControlVarDet}) is only interesting theoretically, because the parameter $p_-^V$ is unknown in practice. We will use (\ref{eq:ControlVarEmp}) for our model selection algorithm. 
\end{itemize}

\subsection{Model Selection}\label{subsection:ModelSelection}
We deduce from Theorem \ref{theo:controlVar} that the risk of the estimator $\Ph_{i|V}$ is bounded in the following way. For all $\delta>1$, for all subsets $V$, 
\begin{equation}\label{eq:RiskBound1}
P\left(\left\|P_{i|G}-\Ph_{i|V}\right\|_{\infty}\leq\left\|P_{i|G}-P_{i|V}\right\|_{\infty}+c_2\sqrt{\frac{\ln(\delta n)}{n\ph^V_-}}\right)\geq 1-\delta^{-1}.
\end{equation}
The risk of $\Ph_{i|V}$ depends on the approximation properties of $V$ through the bias $\left\|P_{i|G}-P_{i|V}\right\|_{\infty}$ that is typically unknown in practice, and on the complexity of $V$, measured here by $\ph^V_-$. The aim of this section is to provide model selection procedures in order to select a subset of $V_M$ that optimizes the bound (\ref{eq:RiskBound1}). In the following, we denote by $\G_n$ a finite collection of subsets of $V_M$, possibly random, and we call optimal or oracle in $\G_n$, any subset $\Gh=\Gh(X_{1:n}(\cup_{V\in \G_n}V))$ in $\G_n$, possibly random, such that,
\begin{equation*}
P\left(\left\|P_{i|G}-\Ph_{i|\Gh}\right\|_{\infty}\leq K\inf_{V\in \G_n}\left\{\left\|P_{i|G}-P_{i|V}\right\|_{\infty}+\sqrt{\frac{\ln(\delta n)}{n\ph^V_-}}\right\}\right)\geq 1-\delta^{-1}.
\end{equation*}
We introduce the following selection rule. Let $N_n$ be an almost sure bound on the cardinality of $|\G_n|$. For all $\delta>1$ and for all $C> c_{2}$, let
\begin{equation}\label{def:estselectinfty}
\Gh(C,\delta,\G_{n})=\argmin_{V\in\G_n}\set{-\left\|\Ph_{i|V}\right\|_{\infty}+C\pen(V)},\;{\rm where}\; \pen(V)\geq \sqrt{\frac{\ln(\delta nN_{n})}{n\ph^V_-}}.
\end{equation}
The following theorem states that $\Gh(C,\delta,\G_n)$ is almost an oracle.
\begin{theorem}\label{theo:GenModSel}
Let $P$ be a probability measure on $\x(G)$ satisfying {\bf H1}. Let $\G_n$ be a finite collection of finite subsets of $G$, possibly random, and let $N_n$ be an almost sure bound on the cardinality of $\G_n$. For all $C>c_{2}$, $\delta>1$, let $\Gh_{\delta}(C)=\Gh(C,\delta,\G_{n})$ be the estimator given by (\ref{def:estselectinfty}). There exists a positive constant $K=K(c_2,C,\kappa_{\min})$ such that,
\begin{equation*}
P\left(\left\|\Ph_{i|\Gh_{\delta}(C)}-P_{i|G}\right\|_{\infty}\leq K\inf_{V\in \G_n}\left\{\left\|P_{i|G}-P_{i|V}\right\|_{\infty}+\pen(V)\right\}\right)\geq 1-\frac1{\delta}.
\end{equation*}
\end{theorem}

\noindent
{\bf Remarks:}
\begin{itemize}
\item Theorem \ref{theo:GenModSel} states that the risk of the estimator selected by the rule (\ref{def:estselectinfty}) is the best among the collection $\G_n$. It is the main result of the paper and we will discuss in what follows several applications.
\item The key idea of the proof is that, by assumption {\bf H1}, we have $\left\|P_{i|G}\right\|_{\infty}-\left\|P_{i|V}\right\|_{\infty}\simeq \left\|P_{i|G}-P_{i|V}\right\|_{\infty}$, hence, our decision rule consists essentially in minimizing the sum of the bias term and the variance term of the risk, and the selected estimator is then an oracle.
\item The constant $c_{2}$ derived in Theorem~\ref{theo:controlVar} is very pessimistic. Hence, Theorem \ref{theo:GenModSel} is more interesting theoretically. In the simulations of Section \ref{sect:simu}, we will calibrate $C$ with the slope algorithm introduced in \cite{BM07} and illustrate the nice properties of the resulting $\Gh_{\delta}(C)$.
\end{itemize}
Let us go back to the ON problem. It is solved thanks to the following corollary.

\begin{coro}\label{coro:ONProblem}
Let $(G,A,P)$ be a random field. Let $V_{M}$ is a finite subset of $G$ with cardinality $M$, let $\delta>1$ and let $\Gamma_{M}(\delta)=\ln(n)(1+\log_{2}(M))+\ln(\delta)$. For all $m$, $e\leq m\leq M$, let $\G_{m,M}=\left\{V\subset V_{M},\; |V|\leq m\right\}.$ For all $V\subset V_{M}$, let $\pen(V)= (n\ph^V_-)^{-1/2}\sqrt{\Gamma_{M}(\delta)}$, let $\Gh_{\delta}(C)=\Gh(C,\delta,\G_{\log_{2(n)},M})$ be the estimator given by (\ref{def:estselectinfty}). With probability larger than $1-\delta^{-1}$, we have
\begin{equation*}
\norm{\Ph_{i|\Gh_{\delta}(C))}-P_{i|G}}_{\infty}\leq K\inf_{V\subset V_{M}}\left\{\left\|P_{i|G}-P_{i|V}\right\|_{\infty}+\sqrt{\frac{\Gamma_{M}(\delta)}{n\ph^V_-}}\right\}.
\end{equation*}
\end{coro}

\noindent
{\bf Remarks:}
\begin{itemize}
\item The complexity of the model selection algorithm for the collection $\G_{m,M}$ is $O(nM^{m})$. This collection is used when a uniform bound $m$ on the cardinalities of the $|G_{i}|$ is known. The complexity is the minimal necessary to recover the interaction graph in this problem \cite{Bresler08}.
\end{itemize}

\subsection{Estimation of the interaction subgraph}
Let $M$ be an integer and let $V_M$ be a finite subset of $G$, with cardinality $M$. For all subsets $V_{M}$ of $G$, let us choose $v^{V}_n(\delta)\geq \sqrt{\ln(\delta n)}(n\ph^V_-)^{-1/2}.$
Let $V$ be a finite subset of $G$, we study in this section the estimators of $G_i$ given by
\begin{align}\label{eq:estGo}
\Gh^{V}_{i}(c)=\left\{j\in V,\;  \omega^{V}_{i,j}(\Ph)>c v^{V}_n(\delta)\right\}.
\end{align}

We introduce the following function.
\begin{equation*}
\Psi(v)=\inf_{V,\;\ph_-^V\geq v^{-2}}\norm{P_{i|V}-P_{i|G}}_{\infty}.
\end{equation*}
$\Psi$ represents the minimal value of the bias term at a given value of the variance term. Our assumption concerns the rate of convergence of $\Psi$ to $0$.

\vspace{0.2cm}

\noindent
{\bf H2($\epsilon_{\Psi}$):} {\it There exist $C_{\Psi}>0$, $\alpha_{\Psi}>0$ such that, for all $K>1$, for all $v>0$,} 
$$P\paren{\Psi(K v)\leq C_{\Psi}K^{-\alpha_{\Psi}}\Psi(v)}\geq 1-\epsilon_{\Psi}.$$

\begin{theorem}\label{theo:IntGraphEst}
Let $(G,A,P)$ be a random field satisfying {\bf H1, H2}. Let $e\leq M$ be an integer, let $V_M$ be a finite subset of $G$ with cardinality $M$. Let $\delta>1$ and let $\Gamma_{M}(\delta)=\ln(n)(1+\log_{2}(M))+\ln(\delta)$. Let $\G_{n}=\set{V\subset V_{M},\;|V|\leq (\log_{2} n)}$. For all $V$ in $\G_{n}$, let 
$$v_n^V(\delta)=\sqrt{\frac{\Gamma_{M}(\delta)}{n\ph_{-}^{V}}}.$$
Let $C\geq c_{2}$, $\pen(V)=v_n^V(\delta)$ and let $\Gh_{\delta}(C)=\Gh(C,\delta,\G_{n})$ be the set selected by the selection rule (\ref{def:estselectinfty}). Let $c>0$ and let $\Gh^{\Gh_{\delta}(C)}_{i}(c)$ be the associated set defined by (\ref{eq:estGo}). Let $K$ be the constant defined in Theorem \ref{theo:GenModSel} and let $c_{\infty}=2\paren{c_{2}+C_{\Psi}^{-1/\alpha_{\Psi}}(2K)^{1-1/\alpha_{\Psi}}}$. We have
\begin{align*}
P\left(\right.&\set{j\in V_M,\;\omega^{G}_{i,j}(P)\geq (c+c_{\infty})v_n^{\Gh}(\delta)}\\
&\left.\subset \Gh^{\Gh_{\delta}(C)}_{i}(c) \subset \set{j\in V_M,\;\omega^{G}_{i,j}(P)\geq (c-c_{\infty})v_n^{\Gh}(\delta)}\right)\geq 1-\delta^{-1}-\epsilon_{\Psi}.
\end{align*}
\end{theorem}

\noindent
{\bf Remark:}
\begin{itemize}
\item When $c>c_{\infty}$, $\Gh^{\Gh_{\delta}(C)}_{i}(c)$ contains exactly the sites that have a pairwise interaction with $i$ of order the risk of an oracle. It provides a partial solution to the INI problem.
\item Theorem \ref{theo:IntGraphEst} requires the extra assumption {\bf H2} compared to Theorem \ref{theo:GenModSel}. Moreover, the theoretical constant $c_{\infty}$ depends on the constants $\kappa_{\min}$, $C_{\Psi}$, $\alpha_{\Psi}$.
\end{itemize}
Let us conclude this section with the two steps algorithm suggested by Theorem \ref{theo:IntGraphEst} to estimate $G_i=\{j\in G,\;\omega^G_{i,j}(P)>0\}$.

\vspace{0.5cm}

\noindent
{\bf Estimation algorithm:}
\begin{itemize}
\item Choose a large subgraph $V_M$ of $G$, typically the $M$ nearest neighbors of $i$ in $G$.
\item {\bf Selection step.} Choose a model $\Gh$, applying the model selection algorithm of Theorem \ref{theo:GenModSel} to the collection of all subgraphs of $V_M$ with cardinality smaller than $\log_2(n)$.
\item {\bf Cutting step.} Cut the edges of $\Gh$ such that $\omega^{\Gh}_{i,j}(\Ph)>c_{\infty}v_n^{\Gh}$.
\end{itemize}

\section{Ising Models}\label{sect:pairGibbsMeasures}
The remaining of the paper is devoted to Ising models. These models are very important in statistical mechanics \citep{Giorgii88} and neuroscience \citep{Schneidman06} where they represent the interactions respectively between particles and neurons. In this section, we prove that Ising models satisfy {\bf H1}, so that all our general results apply in these models. We also define effective algorithms for the ON and INI problems, adapted to this special case.
\subsection{Verification of {\bf H1}.}
Let us recall the definition of Ising models.
\begin{definition}
Let $f: G^2\times A^2\rightarrow \R,\; (i,j,a,b)\mapsto f_{i,j}(a,b)$ be a real valued function. For all $i,j$ in $G$ and all $a$ in $A$, let $\|f_{i,j}^a\|=\max_{b\in A}|f_{i,j}(a,b)|.$ $f$ is said to be a pairwise potential of interaction if, for all $a,b$ in $A$, $f_{i,i}(a,b)=0$ and if
$$r:=\sup_{i\in G}\sup_{a\in A}\sum_{j\in G}\|f_{i,j}^a\|<\infty.$$
In this case, $T=r^{-1}$ is called the temperature parameter of the pairwise potential $f$.
\end{definition} 
\begin{definition}
A probability measure $P$ on $\x(G)$ is called an Ising model with potential $f$ if, for all $x\in \x(G)$,
$$P_{i|G}(x)=\frac{e^{\sum_{j\in G}f_{i,j}(x(i),x(j))}}{\sum_{a\in A}e^{\sum_{j\in G}f_{i,j}(a,x(j))}}=\frac{1}{1+e^{\sum_{j\in G}f_{i,j}(x_i(i),x(j))-f_{i,j}(x(i),x(j))}}.$$
\end{definition}
\noindent The existence of a such a measure is well known \cite{Giorgii88}.
\\

\noindent
{\bf Remark:}
\begin{itemize}
\item The classical Ising model has potential $f$ defined by $f_{ij}(a,b)=J_{ij}ab + H_ia\mathbf{1}_{\{i=j\}}$, $J_{ij}\in \R$, $ H_i \in \R$, for all $a,b \in A$ and $i,j \in G$.
\item One of the fundamental questions studied for this class of models is the description of conditions on potential $f$ that guarantees uniqueness and non-uniqueness of the Ising model. Usually, high temperature implies conditions for the uniqueness of the Ising model and low temperature implies non-uniqueness \cite{Giorgii88}.
\end{itemize}

Let $g_{i,j}(a,b)=f_{i,j}(a,b)-f_{i,j}(-a,b)$, we have then
$$P_{i|G}(x)=\frac{1}{1+e^{-\sum_{j\in G}g_{i,j}(x(i),x(j))}}.$$
It is clear that Ising models satisfy {\bf CA} and {\bf NN} with $\nu=(1+e^{2r})^{-1}$.

\begin{definition}
Let $(G,A,P)$ be an Ising model, with potential $f$. For all $i,j$ in $G$, for all $a$ in $A$, let
$$\omega_{i,j}(f)=\sup_{(a,b)\in A^2}\left\{g_{i,j}(a,b)-g_{i,j}(a,-b)\right\}=\sup_{b\in A}\left\{g_{i,j}(a,b)-g_{i,j}(a,-b)\right\}.$$
\end{definition}

\noindent
Let us first recall some elementary facts about Ising models.
\begin{proposition}\label{prop:elGibbs}
Let $(G,A,P)$ be an Ising model, with potential $f$. For all finite subsets $V$ of $G$, for all $i,j$ in $G$, we have
\begin{enumerate}
\item $p_-^V\geq (1+e^{2r})^{-|V|}$.
\item $\frac{2e^{-2r}}{(1+e^{2r})^2}\omega_{i,j}(f)\leq \omega^{G}_{i,j}(P)\leq \frac{e^{2r}(e^{4r}-1)}{4r(1+e^{-2r})^2}\omega_{i,j}(f)$.
\end{enumerate}
\end{proposition}

\noindent
The following theorem states that all of our general results apply in Ising models. The key ingredient of the proof is the precise control of the bias term (\ref{eq:ControlBias1}).
\begin{theorem}\label{theo:controlBiasGibbs}
Let $(G,A,P)$ be an Ising model, with potential $f$. There exist two positive constants $c_r^*\leq C_r^*$ such that, for all subsets $V$ of $G$,
\begin{equation}\label{eq:ControlBias1}
c_r^*\sum_{j\notin V}\omega_{i,j}(f)\leq \left\|P_{i|G}-P_{i|V}\right\|_{\infty}\leq C_r^*\sum_{j\notin V}\omega_{i,j}(f).
\end{equation}
$P$ satisfies assumption {\bf H1} i.e. there exists a constant $\kappa_{\min}>0$ such that, for all finite subsets $V$ of $G$,
$$\kappa_{\min}\left\|P_{i|G}-P_{i|V}\right\|_{\infty}\leq \left\|P_{i|G}\right\|_{\infty}-\left\|P_{i|V}\right\|_{\infty}.$$
\end{theorem}

\subsection{A special strategy for Ising models} \label{sec:reduction}
The model selection algorithm (\ref{def:estselectinfty}) might be computationally demanding in practice when the collection $\G_n$ is too large. This is the case of the collection $\G_{\log_2(n),M}$ used several times in Section \ref{Section:GeneralResults}, when the values of $M$ and $n$ are large. The purpose of this section is to show that a special strategy, computationally more attractive, can be adopted in Ising models. The idea comes from \cite{Bresler08}. Let us describe the method.

\vspace{0.2cm}

\noindent
\emph{Reduction of the number of sites.} Let $x_1$ be the configuration in $\x(G)$ such that, for all $j$ in $G$, $x_1(j)=1$.

\begin{itemize}

\item[Step 1] Computation of the empirical probabilities. For all $j$ in $V_M$, let 
$$\ph(j)=\Ph(x_1(j)),\; \ph(i,j)=\Ph(x_1(i,j)).$$

\item[Step 2] Reduction step. We keep the $j$ in $V_M$ such that
$$|\ph(i,j)-\ph(i)\ph(j)|>\eta.$$
\end{itemize}

Let also $\eta_{ms}$ be the smallest $\eta>3\sqrt{(2n)^{-1}\ln(6M\delta)}$ such that the number of $j$ kept after Step 2 is smaller than $\kappa\log_{2}(n)$.

We denote by $\Vh(\eta)$ the set of $j$ kept after Step 2. It is clear that the reduction algorithm has a complexity $O(nM)$. Remark that the values $|\ph(i,j)-\ph(i)\ph(j)|$ do not depend on the configuration $x_1$ since the alphabet has only two letters.
\vspace{0.2cm}

\noindent
\emph{Model selection algorithm.} Let $\G=\set{V\subset \Vh(\eta_{ms})}$. 

\begin{itemize}

\item[Step 1] Computation of the conditional probabilities. For all $V$ in $\Vh(\eta_{ms})$, compute $\norm{P_{i|V}}$, and $\pen(V)$.

\item[Step 2] Selection Step. We choose $C>c_2$ and
$$\Gh=\arg\min_{V\in \G}\set{-\norm{P_{i|V}}+C\sqrt{\frac{\ln(n^{\kappa}\delta)}{n\ph_-^{V}}}}.$$
\end{itemize}

It is clear that, if $\mh=|\Vh(\eta_{ms})|\leq \kappa (\log_2(n))$, hence
$$\Nh=|\G|=\sum_{k=0}^{\mh}C_{\mh}^k\leq 2^{\mh}\leq n^{\kappa}.$$
Hence, the complexity of the model selection algorithm is $O(n^{\kappa+1})$. The global complexity of the algorithm is therefore $O(n^{\kappa+1}+nM)$. As a comparison, the model selection algorithm for $\G_n=\G_{\log_2(n),M}$ was $O(nM+n^{\log_2(M)})$.

\subsubsection{Control of the risk of the resulting estimator}

\begin{theorem}\label{theo:BornesRisk}
Let $(G,A,P)$ be an Ising model, with potential $f$. Let 
$$C_{1}=\frac{4r(1+e^{2r})^{3}}{e^{-6r}(e^{4r}-1)}, C_{2}=\frac{4r(1+e^{2r})^{2}}{e^{6r}(e^{4r}-1)}.$$
With probability larger than $1-\delta$ we have that
\begin{align*}
&\set{j\in V_{M},\;\absj{\omega_{i,j}(f)}\geq C_{1}\paren{\eta+3\sqrt{\frac{\ln(6M\delta)}{2n}}}}\\
&\subset\Vh(\eta)\subset\set{j\in V_{M};\; \absj{\omega_{i,j}(f)}\geq C_{2}\paren{\eta-3\sqrt{\frac{\ln(6M\delta)}{2n}}}}.
\end{align*}
Furthermore, let us denote by 
$$V(\delta,M)=\set{j\in V_M;\;|\omega_{i,j}(f)|\leq C_{1}(\eta_{ms}+3\sqrt{(2n)^{-1}\ln (6M\delta)}}.$$
With probability larger than $1-2\delta$, we have,
$$\frac1K\norm{\Ph_{i|\Gh}-P_{i|G}}\leq \sum_{j\in V(\delta,M)}|\omega_{i,j}(f)|+\inf_{V\in\G}\set{\sum_{j\in \Vh(\eta)/V}|\omega_{i,j}(f)|+\sqrt{\frac{\ln(n^{\kappa}\delta)}{n\ph_-^{V}}}}.$$
\end{theorem}

\noindent
{\bf Remarks:}
\begin{itemize}
\item The estimator of the interaction graph has better properties than the one obtained with selection and cutting procedure. The main difference is that there is no term $(\ph_-^{\Gh})^{-1/2}$ in the rate of convergence.
\item The oracle inequality might be a little bit less sharp than the one obtained in (\ref{eq:ModSelGibbs}). This is the price to pay to have a computationally efficient algorithm.
\item Our result holds in the Ising model. However, \cite{Bresler08} used a similar approach in more general random fields with some additional assumptions and obtained good properties for the INI problem.
\end{itemize}

\section{Simulation studies}\label{sect:simu}
\begin{figure}[!ht]
\begin{center} 
 \includegraphics[width=6cm]{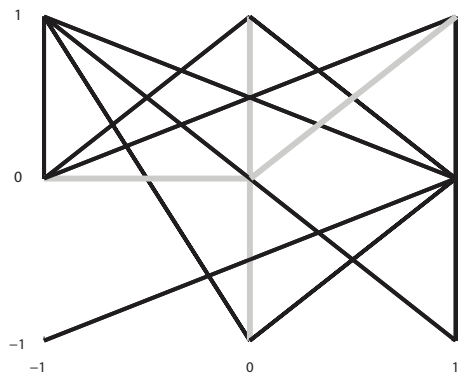}
 \caption{\label{fig:InteractionMatrix_2_weak} \small Representation of the interacting pairs of the Ising model used in the simulation experiments. The edges between sites indicate the interacting pairs. The grey colored edges indicate the sites interacting with site $(0,0)$. } \label{fig:intmatrix2weak}
 \end{center}
\end{figure}

In this section we illustrate results obtained in Sections \ref{Section:GeneralResults} and \ref{sect:pairGibbsMeasures} using simulation experiments and introduce the slope heuristic. All these simulation experiments can be reproduced by a set of MATLAB\textsuperscript{\textregistered} routines that can be downloaded from \href{http://www.princeton.edu/~dtakahas/publications/LT10routines.zip}{www.princeton.edu/$\sim$ dtakahas/publications/LT10routines.zip}. \\
Let $G=\{-1, 0, 1\}\times \{-1, 0, 1\}$. For the sections \ref{sec:var}, \ref{sec:slope}, \ref{sec:oracriskcomp}, \ref{sec:discoveryrateON}, \ref{sec:performanceINI}, \ref{sec:relationshipINION}, and \ref{sec:seleccut}, we consider an Ising model on $A^G$, with pairwise potential given by $f_{ij}(c,d)=J\mathbf{1}_{j\in V_i}cd$ for $i,j \in G$, $c,d \in A$, $J = 0.2$, and $V_i \subset G$. The pair of sites $(i,j)$ where $j \in V_i$ is shown in Figure \ref{fig:intmatrix2weak}.  For all these experiments, $i = (0,0)$. We simulated independent samples of the Ising model with  increasing sample sizes $n=100k$, $k=1, \ldots, 100$. For each sample size we have $N=100$ independent replicas.

\subsection{Variance term of the risk} \label{sec:var}
\begin{figure}[!ht]
\begin{center}
 \includegraphics[width=10cm]{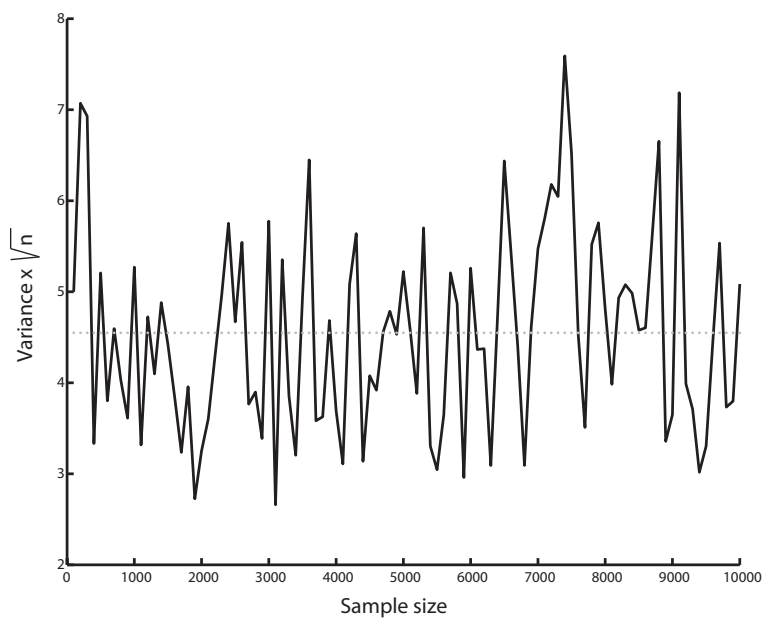}
 \caption{\label{fig:SupNorm1_J02_reduced} \small Plot of the number of samples $n$ against $\sqrt{n}\left\|\hat{P}_{i|V_i}-P_{i|V_i}\right\|_\infty $. The dotted line indicates the linear regression line. Observe that the regression line is essentially parallel to the abscissa.} 
 \end{center}
\end{figure}

In the following experiment we will verify Theorem \ref{theo:controlVar} in a simulation.
For each sample size we computed the normalized variance term $\sqrt{n}\left\|\hat{P}_{i|V_i}-P_{i|V_i}\right\|_\infty $ for $N$ different samples and obtained the average value.  The result is summarized in Figure \ref{fig:SupNorm1_J02_reduced}.

\subsection{Slope heuristic}\label{sec:slope}

The constant $c_2$ derived from Theorem \ref{theo:controlVar} is too pessimistic to be used in practice. The purpose of this section is to present a general method to design this constant. It is based on the slope heuristic, introduced in \cite{BM07} and proved in several other frameworks in \cite{ AM08, Le09}. We refer also to \cite{ Michel10} for a large discussion on the practical use of this method. In order to describe it, let us introduce, for all $V$ in $\G_{m,M}$, a quantity $\Delta_V$, possibly random, measuring the complexity of the model $V$. The heuristic states the following facts.

\begin{enumerate}
\item There exists a positive constant $C_{\min}$ such that when $C < C_{\min}$, the complexity of the model selected by the rule (\ref{def:estselectinfty}) is as large as possible. 
\item When $C $ is slightly larger than $C_{\min}$  the complexity of the selected model is much smaller. 
\item When $C = 2C_{\min}$ then the risk of the selected model is asymptotically the one of an oracle. 
\end{enumerate}

\noindent
The heuristic yields the following algorithm, defined for all complexity measures $\Delta_V$.

\begin{enumerate}
\item For all $C>0$, compute $\Delta_{\hat{G}(C)}$, the complexity of the model selected by the rule (\ref{def:estselectinfty}). 
\item Choose $\tilde{C}_{\min}$ such that $\Delta_{\hat{G}(C)}$ is very large for $C< \tilde{C}_{\min}$ and much smaller for $C> \tilde{C}_{\min}$.
\item Select the final $\hat{G}=\hat{G}(2\tilde{C}_{\min})$. 
\end{enumerate}

The algorithm is based on the idea that $\tilde{C}_{\min}\simeq C_{\min}$ and therefore that the final $\hat{G}$, selected by $2\tilde{C}_{\min}\Delta_V$ is an oracle by the third point of the slope heuristic. The actual efficiency of this approach depends highly on the choice of the complexity measure $\Delta_V$ and on the practical way to choose $\tilde{C}_{\min}$ in step 2 of the algorithm. 
We illustrate the dependence on $\Delta_{V}$ in the following experiences. 

$\Delta_V$ is either the cardinality of $V$ (the dimension) or the variance estimator $C(n\hat{p}_-^{V})^{-1/2}$. $\tilde{C}_{\min}$ is selected with the maximum jump criteria \cite{ Michel10}: fix an increasing sequence of positive numbers $C_0, \ldots, C_t$ and define
\begin{equation*}
k = \argmax_i \left\{ \Delta_{\hat{G}(C_{i})}-\Delta_{\hat{G}(C_{i-1})} \right\},\;{\rm and}\; \tilde{C}_{\min} = C_k. 
\end{equation*}
If the maximum is achieved in more than one value, take the biggest of such $k$.

\vspace{0.2cm}

\noindent
{\bf Remark:} The calculation of $\tilde{C}_{\min}$ does not yield a significant increase of computational time compared to the evaluation of the model selection criteria for one fixed constant $C$. The only additional cost is due to the fact that one has to keep in the computer memory the conditional probabilities that must be computed only once.

\subsection{Oracle risk compared to the risk of the estimated model} \label{sec:oracriskcomp}
\begin{figure}[!ht]
\begin{center}
 \includegraphics[width=10cm]{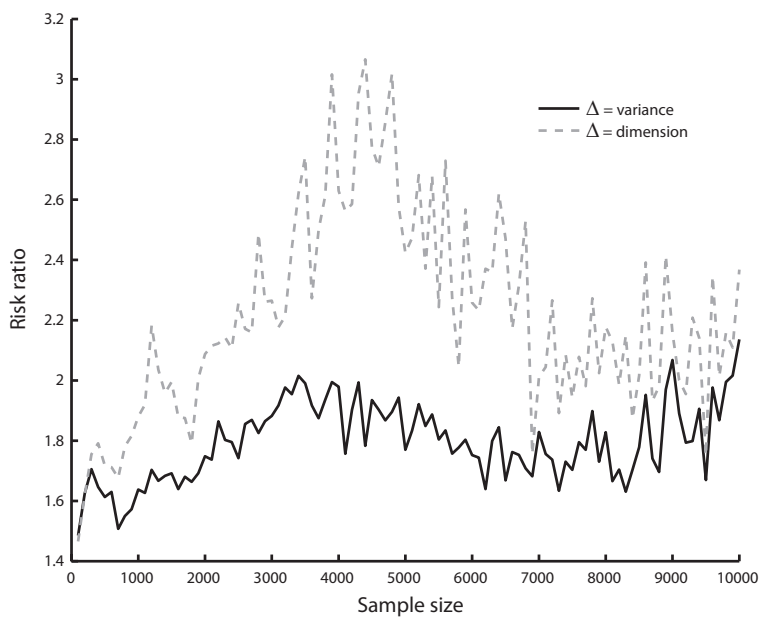}
 \caption{\label{fig: RiskRatio_weak_3X3} \small Plot of the number of samples $n$ against the average of ratio (\ref{eq:riskratio}). Observe that the risk ratio remains bounded for both the variance (solid black) and the dimension (dashed grey) as the measure of complexity. } 
 \end{center}
\end{figure}

One way to verify the performance of the slope heuristic proposed in previous section is to compute the ratio
\begin{equation} 
 \frac{\left\|\hat{P}_{i|\Gh(2\tilde{C}_{\min})}-P_{i|G}\right\|_\infty }{\inf_{V \subset G}\left\|\hat{P}_{i|V}-P_{i|G}\right\|_\infty}. \label{eq:riskratio}
 \end{equation}
With a reasonable procedure, we expect that the above quantity remains bounded. We applied the model selection procedure (\ref{def:estselectinfty}) with slope heuristic discussed above for the set $\{V\subset G\setminus \{i\}: |V| \leq 8\}$.
For each sample size we computed the ratio (\ref{eq:riskratio}) for 100 different samples and we obtained the average. The result is summarized in Figure \ref{fig: RiskRatio_weak_3X3}.

\subsection{Discovery rate of the model selection procedure for ON problem} \label{sec:discoveryrateON}
\begin{figure}[!ht]
\begin{center}
 \includegraphics[width=10cm]{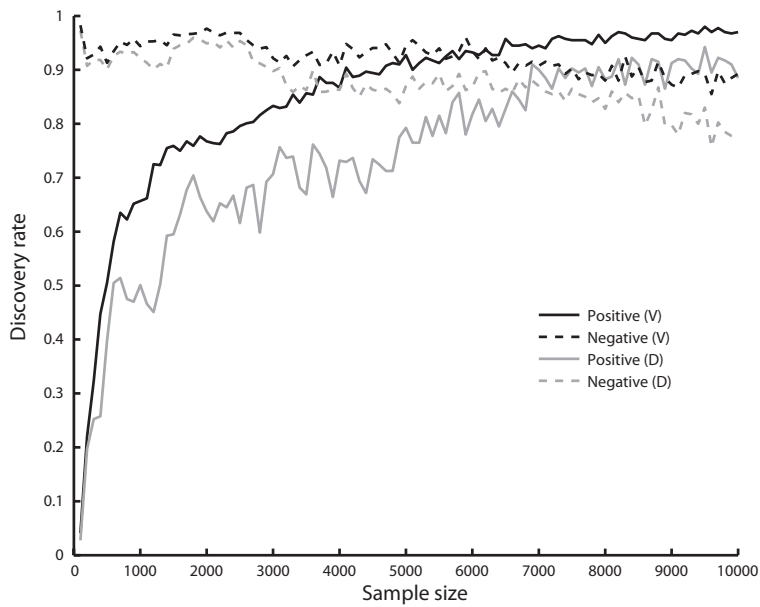}
 \caption{\label{fig:Figure4_LT10} \small Plot of positive and negative discovery rates with respect to the oracle against the sample size $n$. In solid/dashed black lines are represented the positive/negative discovery rates using the variance (V) as the complexity measure and in solid black/grey lines the positive/negative discovery rates using the dimension (D). Observe that the variance gives a better positive and negative discovery rates with respect to oracle when compared to the dimension. } 
 \end{center}
\end{figure}

Another way to measure the performance of our model selection procedure is to compute the positive discovery rate
\begin{equation} \label{eq:positiveoracle}
\E\left[\frac{|\Gh(2\tilde{C}_{\min}) \cap \Gh_{oracle}|}{|\Gh_{oracle}|} \right]
\end{equation}
and the negative discovery rate
\begin{equation} \label{eq:negativeoracle}
\E\left[\frac{|G\setminus \left(\Gh(2\tilde{C}_{\min}) \cup \Gh_{oracle}\right)|}{|G\setminus \Gh_{oracle}|} \right].
\end{equation}
with respect to the oracle $\Gh_{oracle}$.\\
We estimated (\ref{eq:positiveoracle}) and (\ref{eq:negativeoracle}) and the result is summurized in Figure \ref{fig:Figure4_LT10}.

\subsection{Performance of the model selection procedure for INI problem} \label{sec:performanceINI}
\begin{figure}[!ht]
\begin{center}
 \includegraphics[width=8cm]{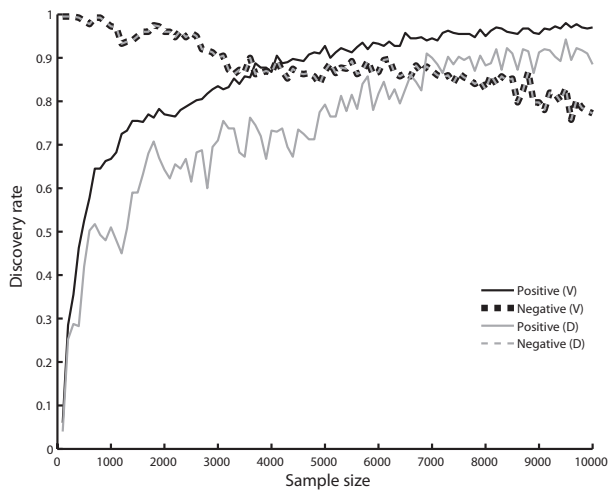}
 \caption{\label{fig:Figure5_LT10} \small Plot of positive and negative discovery rates with respect to $V_i$ against the sample size $n$. In solid/dashed black lines are represented the positive/negative discovery rates using the variance (V) as the complexity measure and in solid black/grey lines the positive/negative discovery rates using the dimension (D). Observe that the variance gives higher positive discovery rates than the dimension as the measure of complexity although the negative discovery rates are the same.} 
 \end{center}
\end{figure}

A natural question is how well the proposed model selection procedure behaves for the INI problem. Observe that the model selection procedure was designed to solve the ON problem and in principle does not necessary work for the INI problem. To investigate this question for each sample size we estimated the positive discovery rate
\begin{equation*}
\E\left[\frac{|\Gh(2\tilde{C}_{\min}) \cap V_i |}{|V_i|} \right]
\end{equation*}
and the negative discovery rate
\begin{equation*}
\E\left[\frac{|G\setminus \left(\Gh(2\tilde{C}_{\min}) \cup V_i\right)|}{|G\setminus V_i|} \right],
\end{equation*}
with respect to the interaction neighborhood $V_i$.
The result is summurized in Figure \ref{fig:Figure5_LT10}.

\subsection{Relationship between the INI and ON problems} \label{sec:relationshipINION}

\begin{figure}[!ht]
\begin{center}
 \includegraphics[width=8cm]{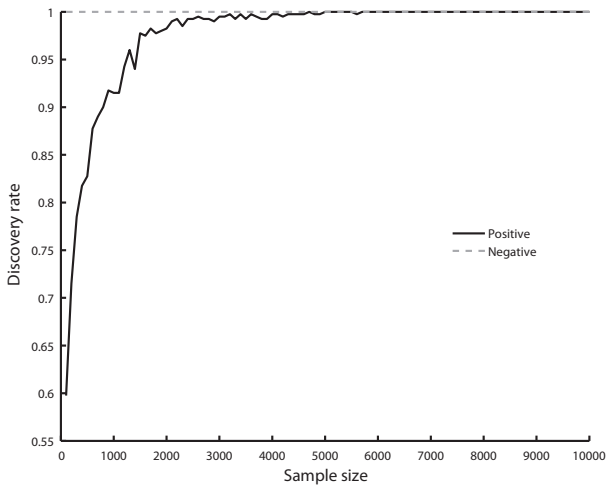}
 \caption{\label{fig:Figure6_LT10} \small Plot of positive and negative discovery rates of the oracle with respect to $V_i$ against the sample size $n$. The solid black line represents the results for positive discovery rates and the dashed grey line represents the results for the negative discovery rates. Observe that in this example the oracle $\Gh_{oracle}$ matches the interaction neighborhood $V_i$ quite fast. Also observe that in this example the oracle never included interactions not contained in $V_i$.} 
 \end{center}
\end{figure}

Another interesting question is to understand what is the relationship between the INI and ON problems. Useful quantities for this are the positive discovery rate
\begin{equation}\label{eq:positiveINIINO}
\E\left[\frac{|\Gh_{oracle} \cap V_i|}{|V_i|} \right]
\end{equation}
and the negative discovery rate
\begin{equation}\label{eq:negativeINIINO}
\E\left[\frac{|G\setminus \left(\Gh_{oracle} \cup V_i\right)|}{|G\setminus V_i|} \right].
\end{equation}\\
We estimated these quantities and the results are summarized in Figure \ref{fig:Figure6_LT10}.

\subsection{Select and cut procedure} \label{sec:seleccut}

Here we will show the usefulness of the two-step procedure introduced in Theorem~\ref{theo:IntGraphEst} by an example. We consider the same independent samples used in previous experiments. We also consider $i = (0,0)$ and sample sizes $n=100k$, $k=1, \ldots, 100$ with $100$ independent replicas for each sample size.

\begin{figure}[!ht]
\begin{center}
 \includegraphics[width=7cm]{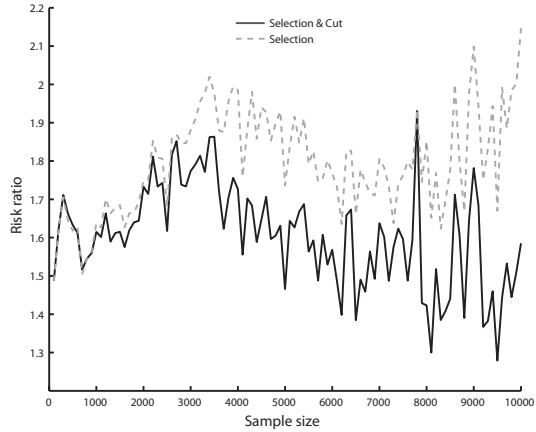}
 \caption{\label{fig:Figure7_LT10} \small Plot of the number of samples $n$ against the average of risk ratio (\ref{eq:riskratiocut}) and (\ref{eq:riskratio}). In solid black is represented the risk ratio for the two-step procedure and in dashed grey the risk ratio for the model selection procedure alone. Observe that the risk ratio of the two-step procedure remains closer to one when compared to the model selection alone.} 
 \end{center}
\end{figure}

Let $\Gh(2\tilde{C}_{\min})$ be the subset of $G$ chosen by first applying the model selection procedure for the set $\{V \subset G\setminus \{i\}: |V| \leq 8\}$. To choose the constant in the model selection procedure, we used the slope heuristic with variance as the complexity measure.  
Let $\Gh(SC)$ be the subset of $G$ obtained by applying to the subset $\Gh(2\tilde{C}_{\min})$ the cutting procedure with $cv^{V}_n=0.3(n\ph_-^V)^{-1}$.  We first computed the average of the risk ratio
\begin{equation} 
 \frac{\left\|\hat{P}_{i|\Gh(SC)}-P_{i|G}\right\|_\infty }{\inf_{V \subset G}\left\|\hat{P}_{i|V}-P_{i|G}\right\|_\infty}. \label{eq:riskratiocut}
 \end{equation}
 for each sample size and compared them with the average of risk ratio (\ref{eq:riskratio}). The results are summarized in Figure \ref{fig:Figure7_LT10}.

\begin{figure}[!ht]
  \begin{center}
 \includegraphics[width=9cm]{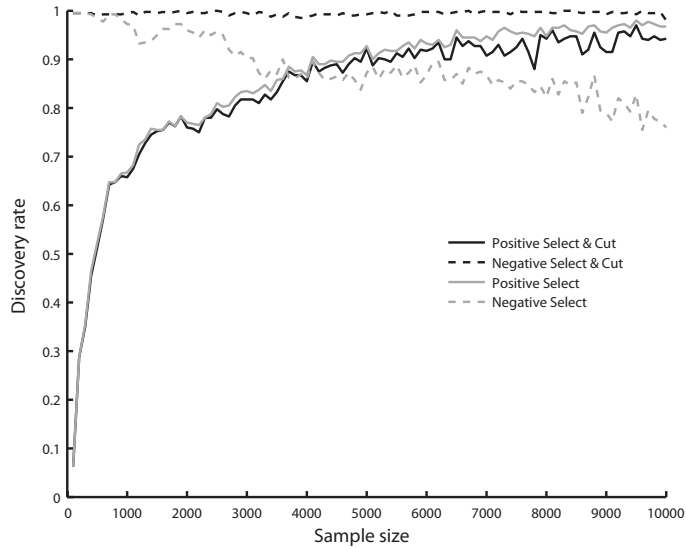}
 \caption{\label{fig:Figure8_LT10} \small Plot of positive and negative discovery rates of $\Gh(SC)$ and $\Gh(2\tilde{C}_{\min})$ with respect to $V_i$ against the sample size $n$. The black solid/dashed lines represent the positive/negative discovery rates of the two-step procedure. The grey solid/dashed lines represent the positive/negative discovery rates of the model selection procedure alone. Observe that the two-step procedure has almost perfect negative discovery rates with incresing positive discovery rates.} 
 \end{center}
\end{figure}

We also computed the positive and negative discovery rates of $\Gh(SC)$ and $\Gh(2\tilde{C}_{\min})$ with respect to $V_i$. The results are presented in Figure \ref{fig:Figure8_LT10}.

\subsection{Computationally efficient algorithm} \label{sec:compeff}
\begin{figure}[!ht]
\begin{center}
 \includegraphics[width=8cm]{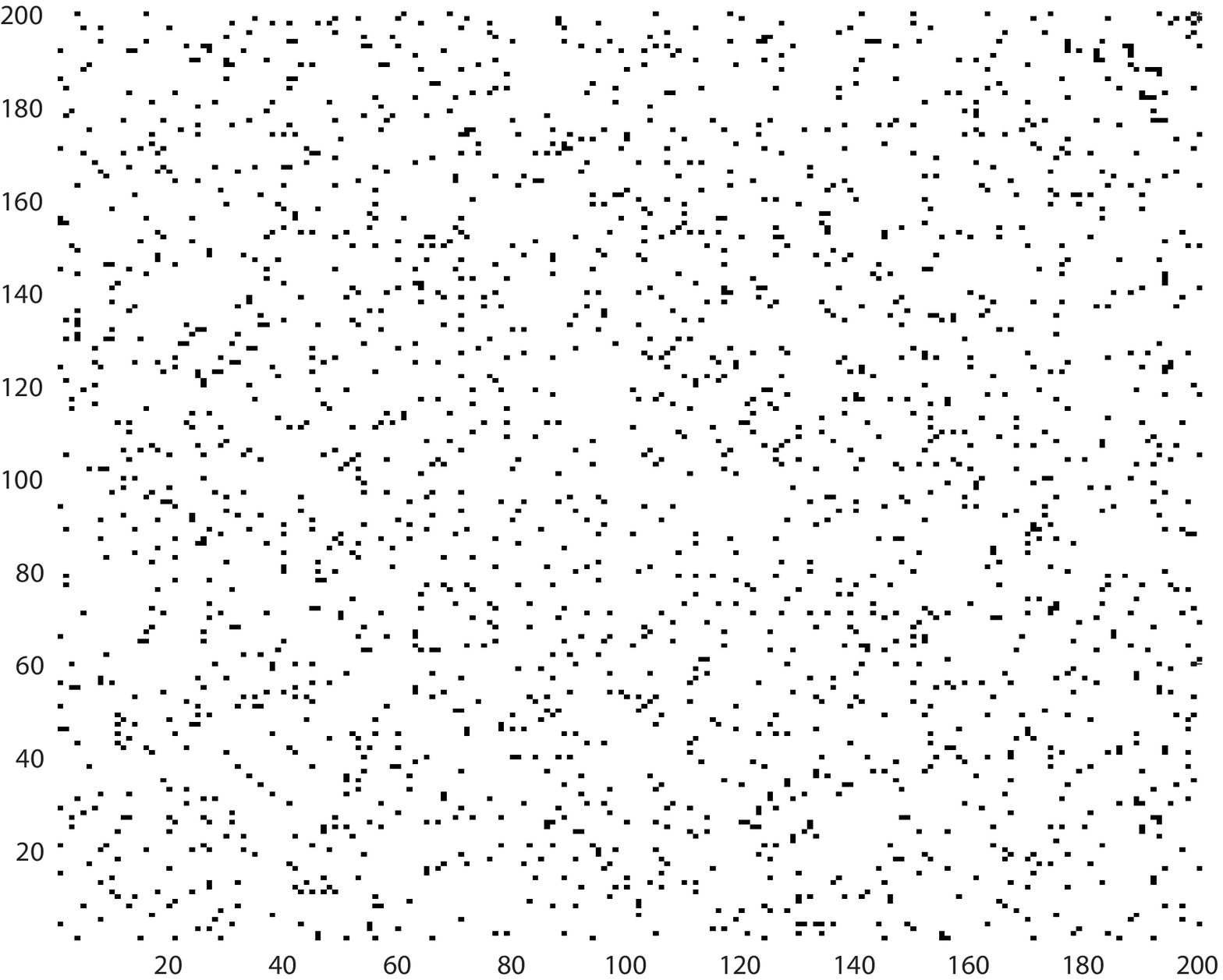}
 \caption{\label{fig:InteractionMatrix_Big} \small Representation of the interacting sites in the Ising model described in \ref{sec:compeff}. The positions $(i,j)$ of the dots indicate the pair of sites $(i,j)$ for which $j\in V_i$.} 
 \end{center}
\end{figure}

In this section we will illustrate the performance of the strategy introduced in Section \ref{sec:reduction} on the Ising model on $A^G$, where $G=\{1, \ldots, 200\}$, with pairwise potential $f_{ij}(c,d)= |J_{ij}|\mathbf{1}_{j\in V_i}cd$ for $i,j \in G$, $c,d \in A$, $V_i \subset G$, and $J_{ij}$ independently generated from a Gaussian distribution with $\E[J_{ij}] = 0$ and $\E[J_{ij}^2] = 4$. The pairs of sites $(i,j)$ with $j \in V_i$ are represented in Figure \ref{fig:InteractionMatrix_Big}.  

\begin{figure}[!ht]
\begin{center}
 \includegraphics[width=10cm]{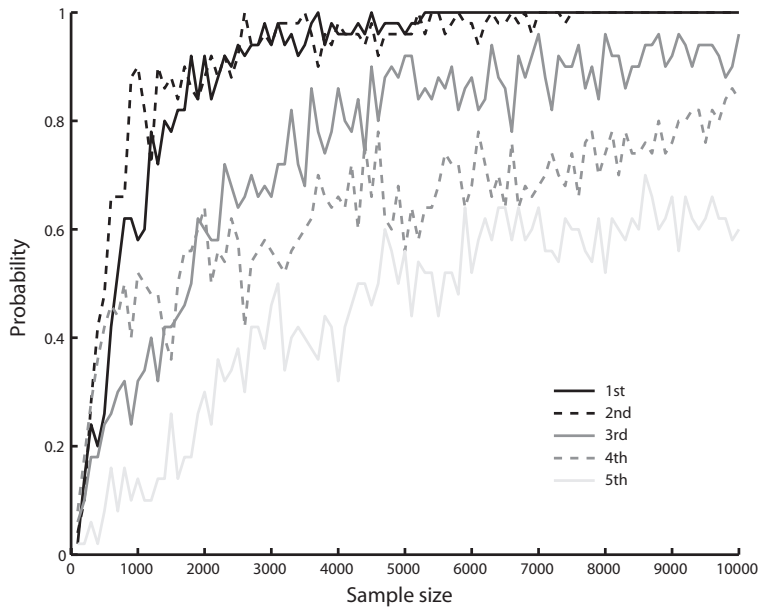}
 \caption{\label{fig:Figure9_LT10} \small Plot of the number of samples $n$ against the probability that  $\Gh_{\mbox{efficient}}$ includes the largest (solid black), and second (dashed black), third (solid gray), fourth (dashed gray), fifth (solid light gray) largest interaction potentials.  Observe that the model selection procedure includes the sites with larger interaction potentials more often.} 
 \end{center}
\end{figure}

For this experiment $i = 1$ and $|V_i| = 16$. We simulated independent samples of the Ising model with  increasing sample sizes $n=100k$, $k=1, \ldots, 100$. For each sample size we have $N=50$ independent replicas. In this example, it is not practical to compute all candidates in collection  $\G_{8,200}$ whereas the algorithm introduced in Section \ref{sec:reduction} is very efficient. We illustrate its performance in the case where the number of sites $j$ kept after Step 2 of the reduction step in Section \ref{sec:reduction} is 10. We denote the model chosen by this algorithm by $\Gh_{\mbox{efficient}}$ We estimated the probability that the selected model $\Gh_{\mbox{efficient}}$ recover the largest, and second, third, fourth, fifth largest interaction potentials. Formally, let $\mathcal{J}_1=\max\{|J_{ij}|\mathbf{1}_{j\in V_i}: i,j \in G\}$ and $\mathcal{J}_k=\max\{|J_{ij}|\mathbf{1}_{j\in V_i}: i,j \in G \setminus \mathcal{J}_{k-1}\}$, for $k = 1, \ldots, 5$. We estimated
  \begin{equation}
  P( \Gh_{\mbox{efficient}} \ni \mathcal{J}_k),
  \end{equation}
  for $k =1, \ldots, 5$. The result of the simulation is presented in Figure \ref{fig:Figure9_LT10}.
  By Monte Carlo simulation using a sample size of $100\,000$ we concluded that the considered Ising model at site $i = 1$ does not satisfy the incoherence condition in \cite{PW10}.

\section{Discussion}\label{Section:Discussion}

We introduced a model selection procedure for interaction neighborhood estimation in partially observed random fields. We prove that the proposed rule satisfies an oracle inequality. The results hold under general assumptions, which for instance, are satisfied by a generalized form of the Ising model. \\
Our model selection approach differs from other works \cite{Bento09, Bresler08, Csiszar06, GOT10, PW10} where only the INI problem is considered and more restrictive conditions are assumed. In particular,  \cite{Bento09, Bresler08, PW10} consider the INI problem for finite random fields and assume that all the interacting sites are observed. This assumption is quite strong from practical point of view, \textit{e.g.} in neuroscience, where the experimenter never has access to the whole set of neurons. Our result holds for partially observed random fields without any restriction on the range of the interactions. \\
Csiszar and Talata \cite{Csiszar06} consider a BIC like consistent model selection procedure for homogeneous random fields in $\Z^d$. As they consider the INI problem using only one realization of the random field, their result are not immediately comparable with ours, but it is interesting to note that they assume that the range of interaction is finite. \\
In \cite{GOT10}, it is considered the INI problem for infinite range Ising models in $\Z^d$. The main restriction in this last work is that it is assumed that the interactions between the sites are weak (``low tempterature")  and that a subset of the observed sites of size $O(\log(n))$, where $n$ is the sample size, must be fixed to apply the proposed procedure. Our procedure has no restriction on the strength of interaction and can be applied for example for low temperature Ising models, provided that the samples come from the same phase. Moreover, the model selection procedure can be applied in high dimension situation when the subset of observed sites is of size $O(n^{\alpha}), \; \alpha>0$.\\ 
We also introduced a two-step procedure in which the model selection rule gives us a small set of candidate sites and a cutting procedure removes from this set the irrelevant interactions. This two step procedure can be understood as a combination of a model selection and a statistical test procedure in spirit of \cite{Zhou10}.\\ 
Our first simulation experiment shows that  the concentration bound for the variance term of the risk in Theorem \ref{theo:controlVar} is sharp.
We propose a slope heuristic with maximal jump criteria using the variance or the dimension as a measure of complexity to choose a good constant in the model selection procedure. In our simulation experiment, we measured the performance of the slope heuristic for ON and INI problems. We observed that the variance had a better behavior as a complexity measure than the dimension because 1) the risk ratio was always smaller for the variance compared to the dimension as the measure of complexity, although both risk ratios remained bounded, 2) the estimated positive and negative discovery rates with respect to the oracle were always higher for the variance compared to the dimension as the measure of complexity, 3) also the estimated positive and negative discovery rates with respect to the interaction neighborhood were always higher for the variance compared to the dimension as the measure of complexity. \\
Although at this point the variance seems to be a better choice for the complexity measure, a more comprehensive study must be carried to obtain a definitive conclusion and  we recommend to consider both measures of complexity in practice.
We addressed also in the simulation experiments the relationship between the ON and INI problems and observed that for sufficiently large sample size, both coincide. The two-step procedure introduced in this article was applied in a example where it clearly enhances the performance of the model selection procedure for both the INI and ON problems. Recently, multistep statistical procedures are gaining attention \cite{Zhou10} although only few rigorous results exist. Our result for the two-step procedure is a contribution for this growing field.   
The main drawback of the proposed model selection procedure is its high computational cost which becomes prohibitive when a large number of sites are observed. We introduced a computationally efficient way to overcome this difficulty in the case of Ising model. The new procedure drastically reduces the set of models for which the model selection procedure must be applied, but still keeping the main interacting sites and a good oracle property. In the simulation experiment we show that the proposed algorithm has a good performance even when the number of the observed sites is as big as 200. It must be remarked that the Ising model considered for this experiment does not satisfy the incoherence condition \cite{PW10} and therefore other computationally efficient algorithms as $\ell_1$-penalizations are not guaranteed to be consistent.
Finally, we provide a set of MATLAB\textsuperscript{\textregistered} routines that can be used  to reproduce our experimental results and to carry further simulation and applied studies.

\subsection*{Acknowledgements}
We would like to thank Antonio Galves and Roberto Imbuzeiro for many discussions during the period in which this paper was written.

\section{Proofs}\label{section:Proofs}

\subsection{Proof of Theorem \ref{theo:controlVar}:}
For all $x$ such that $P(x(V/\{i\}))=0$, we have $\Ph(x(V/\{i\}))=0$, thus $P_{i|V}(x)=\Ph_{i|V}(x)$. Hence, we can only consider the configurations $x$ such that $P(x(V/\{i\}))>0$.
Let us first provide some inequalities about conditional probabilities. 
\begin{lemma}\label{lem:condprob}
Let $x \in \x(G)$, let $V$ be a finite subset of $G$ and let $Q,R$ be two probability measures on $\x(V)$ such that $R(x(V/\{i\}))>0$.
\begin{align*}
Q_{i|V}(x)&-R_{i|V}(x)\\
&= \frac{Q(x(V))-R(x(V))+Q_{i|V}(x)\left(R(x(V/\{i\}))-Q(x(V/\{i\}))\right)}{R(x(V/\{i\}))}.
\end{align*}
\begin{align*}
|Q_{i|V}(x)-R_{i|V}(x)|\leq 3\sup_{x\in\x(G),\;R(x(V/\{i\}))\neq 0}\frac{\left|Q(x(V))-R(x(V))\right|}{R(x(V/\{i\})}.
\end{align*}
\end{lemma}

\noindent
{\bf Remark:} In particular, we deduce from this lemma that 
$$\left\|\Ph_{i|V}-P_{i|V}\right\|_{\infty}\leq 3\sup_{x\in\x(G),\;P(x(V/\{i\})\neq 0}\frac{\left|\Ph(x(V))-P(x(V))\right|}{P(x(V/\{i\}))}.$$
The first inequality follows from the fact that $R_{i|V}(x)=R(x(V))/R(x(V/\{i\}))$ and $Q(x(V))=Q_{i|V}(x)Q(x(V/\{i\}))$. The second one is consequence of the first one and the fact that 
$$|R(x(V/\{i\}))-Q(x(V/\{i\}))|\leq |R(x(V))-Q(x(V))|+|R(x_i(V))-Q(x_i(V))|.$$
The proof of (\ref{eq:ControlVarDet}) is concluded thanks to the following Lemma.
\begin{lemma}\label{lem:pn-p/p}
Let $P$ be a probability measure on $\x(G)$ and let $V$ be a finite subset of $G$. Let $\x'(V)=\{x\in \x(G),\;P(x(V/\{i\}))\neq 0\}$, $p^V_-=\inf_{x\in \x'(G)}P(x(V/\{i\}))$. For all $\delta>1$, with probability larger than $1-\delta^{-1}$, we have
\begin{align*}
\sup_{x\in\x'(G)}\frac{\left|\Ph(x(V))-P(x(V))\right|}{P(x(V/\{i\}))}\leq 64\sqrt{2}\sqrt{\frac{\ln(16\delta/p_-^V)}{np^V_-}}+2048\frac{\ln(16\delta/p^V_-)}{np^V_-}.
\end{align*}
\end{lemma}

\noindent
{\bf Conclusion of the proof of (\ref{eq:ControlVarDet})}. We deduce from Lemmas \ref{lem:condprob} and \ref{lem:pn-p/p} that, with probability larger than $1-\delta^{-1}$,
$$\left\|\Ph_{i|V}-P_{i|V}\right\|_{\infty}\leq 192\sqrt{2}\sqrt{\frac{\ln(16\delta/p_-^V)}{np^V_-}}+6144\frac{\ln(16\delta/p^V_-)}{np^V_-}.$$
As this result is trivial when $\frac{\ln(16\delta/p^V_-)}{np^V_-}>1$, we can always assume that $\frac{\ln(16\delta/p^V_-)}{np^V_-}\leq 1$, hence that $\frac{\ln(16\delta/p^V_-)}{np^V_-}\leq \sqrt{\frac{\ln(16\delta/p^V_-)}{np^V_-}}$, thus, with probability larger than $1-\delta^{-1}$, for $c_1=6144+192\sqrt{2}$,
$$\left\|\Ph_{i|V}-P_{i|V}\right\|_{\infty}\leq c_1\sqrt{\frac{\ln(16\delta/p_-^V)}{np^V_-}}.$$
{\bf Proof of Lemma \ref{lem:pn-p/p}:} We apply Bousquet's version of Talagrand's inequality to the class of functions  $\F=\{(P(x(V/\{i\})))^{-1}1_{x(V)}\}$. This inequality is recalled in Appendix. We have $v^{2}\leq (p_{-}^{V})^{-1}$, $b\leq (p_{-}^{V})^{-1}$, hence, for all $\delta>1$, with probability larger than $1-\delta^{-1}$,
\begin{align}\label{prop:concsup}
\sup_{x\in\x'(G)}&\frac{\left|\Ph(x(V))-P(x(V))\right|}{P(x(V/\{i\}))}\nonumber\\
&\leq 2\E\left(\sup_{x\in\x'(G)}\frac{\left|\Ph(x(V))-P(x(V))\right|}{P(x(V/\{i\}))}\right)+\sqrt{\frac{2\ln(\delta)}{np^V_-}}+2\frac{\ln(\delta)}{np^V_-}.
\end{align}
We apply Lemma~\ref{lemma:EsperanceSupVdisjoints} with $A_x=x(V)$, $x\in \x(G)$, $\alpha_x=[P(x(V/\{i\}))]^{-1}$. We have 
$$\alpha^*=\sup_{x\in\x'(G)} [P(x(V/\{i\}))]^{-1}=\frac1{p^V_-},\;p^*=\sup_{x\in\x'(G)}[P(x(V/\{i\}))]^{-2}P(x(V))\leq\frac1{p^V_-}.$$
Hence,
\begin{equation}\label{prop:Esup}
\E\left(\sup_{x\in\x'(G)}\frac{\left|\Ph(x(V))-P(x(V))\right|}{P(x(V/\{i\}))}\right)\leq\frac{32\sqrt{2}}{\sqrt{np^V_-}}\sqrt{\ln\left(\frac{16}{p_-^V}\right)}+\frac{1024}{np^V_-}\ln\left(\frac{16}{p^V_-}\right).
\end{equation}
Lemma \ref{lem:pn-p/p} is then obtained with (\ref{prop:concsup}) and (\ref{prop:Esup}).\\
Let us now turn to the proof of (\ref{eq:ControlVarEmp}). Let $V$ be a finite subspace of $S$. As (\ref{eq:ControlVarEmp}) holds when $\ph_{-}^{V}=n^{-1}$, it remains to prove (\ref{eq:ControlVarEmp}) when, for all $x$ in $\x(V)$, $\Ph(\x(V))>0$. This is done by the following Proposition.
\begin{proposition}\label{prop:ConcInf}
Let $P$ be a probability measure on $\x(G)$, let $V$ be a finite subset of $G$. Let $\x_{n}=\set{x\in \x(G),\;\Ph(x(V/\{i\})\neq 0}$, $\ph_-^V=\inf_{x\in \x_{n}}\Ph(x(V))$. There exists an absolut constant $c_2\leq 400$ such that, for all $\delta>1$,
\begin{equation}\label{eq:controlforallx}
P\left(\exists x\in\x_{n},\;|\Ph_{i|V}(x)-P_{i|V}(x)|> c_2\sqrt{\frac{\ln (\delta n)}{n\Ph(x(V))}}\right)\leq \frac1{\delta}.
\end{equation}
In particular,
\begin{equation*}
P\left(\sup_{x\in\x_{n}}|\Ph_{i|V}(x)-P_{i|V}(x)|> c_2\sqrt{\frac{\ln (\delta n)}{n\ph_{-}^{V}}}\right)\leq \frac1{\delta}.
\end{equation*}
\end{proposition}

\noindent
{\bf Proof of Proposition~\ref{prop:ConcInf}.} Let $n\geq 2$, $\delta>1$, $c_{2}=400$ and let us first remark that we only have to prove (\ref{eq:controlforallx}) on the subset $\x_{n}'\subset \x_{n}$ of all the $x$ in $\x_{n}$ such that $\Ph(x(V))\geq c_{2}^{2}\ln (\delta n)n^{-1}$.
Let $x$ in $\x_{n}'$, then we also have $P(x(V/\{i\})\neq 0$. From Lemma~\ref{lem:condprob}, we have
\begin{align*}
|\Ph_{i|V}(x)-P_{i|V}(x)|\leq \frac{\absj{\Ph(x(V))-P(x(V))}+\absj{\Ph(x(V/\{i\}))-P(x(V/\{i\}))}}{P(x(V/\{i\}))}.
\end{align*}
From Lemma~\ref{lem:condprob}, we also have
\begin{align*}
|\Ph_{i|V}(x)-P_{i|V}(x)|\leq \frac{\absj{\Ph(x(V))-P(x(V))}+\absj{\Ph(x(V/\{i\}))-P(x(V/\{i\}))}}{\Ph(x(V/\{i\}))\vee c_{2}^{2}\ln (\delta n)n^{-1}}.
\end{align*}
We deduce that
\begin{align*}
|\Ph_{i|V}(x)-P_{i|V}(x)|\leq \frac{\absj{\Ph(x(V))-P(x(V))}+\absj{\Ph(x(V/\{i\}))-P(x(V/\{i\}))}}{P(x(V/\{i\}))\vee \Ph(x(V/\{i\}))\vee c_{2}^{2}\ln (\delta n)n^{-1}}.
\end{align*}
Hence, using the elementary inequality $a\vee b\geq \sqrt{ab}$ with $a=\Ph(x(V/\{i\}))$, $b=P(x(V/\{i\}))\vee c_{2}^{2}\ln (\delta n)n^{-1}$, we deduce that
\begin{align*}
|\Ph_{i|V}(x)-P_{i|V}(x)|\leq \frac{\absj{\Ph(x(V))-P(x(V))}+\absj{\Ph(x(V/\{i\}))-P(x(V/\{i\}))}}{\sqrt{\Ph(x(V/\{i\}))\paren{P(x(V/\{i\}))\vee  c_{2}^{2}\ln (\delta n)n^{-1}}}}.
\end{align*}
We have obtain that, for all $x$ in $\x_{n}'$,
\begin{align*}
&\sqrt{\Ph(x(V/\{i\}))}|\Ph_{i|V}(x)-P_{i|V}(x)|\\
&\leq \frac{\absj{\Ph(x(V))-P(x(V))}+\absj{\Ph(x(V/\{i\}))-P(x(V/\{i\}))}}{\sqrt{P(x(V/\{i\}))\vee  c_{2}^{2}\ln (\delta n)n^{-1}}}\\
&\leq 3\sup_{x\in\x_{n}'}\frac{\absj{\Ph(x(V))-P(x(V))}}{\sqrt{P(x(V/\{i\}))\vee  \ln (\delta n)n^{-1}}}\leq 3\sup_{x\in\x(G)}\frac{\absj{\Ph(x(V))-P(x(V))}}{\sqrt{P(x(V/\{i\}))\vee  c_{2}^{2}\ln (\delta n)n^{-1}}}.
\end{align*}
We apply Bousquet's version of Talagrand's inequality to the class of functions 
$$\F=\set{f=\paren{P(x(V/\{i\}))\vee  c_{2}^{2} \ln (\delta n)n^{-1}}^{-1/2}1_{x(V)},\;x\in \x(G)}$$
We have
$$v^{2}=\sup_{f\in \F}\Var (f(X_{1}))\leq 1,\;b=\sup_{f\in\F}\norm{f}_{\infty}\leq \sqrt{ c_{2}^{-2}(\ln (\delta n))^{-1}n}.$$
Hence, for all $\epsilon>0$, with probability larger than $1-\delta^{-1}$, we have
\begin{align*}
\sup_{x\in\x'_{n}}&\sqrt{\Ph(x(V/\{i\}))}|\Ph_{i|V}(x)-P_{i|V}(x)| \\
&\leq3(1+\epsilon)\E\paren{\sup_{f\in\F}\absj{(P_{n}-P)f}}+3\sqrt{\frac{2\ln(\delta)}n}+\paren{1+\frac3{\epsilon}}\frac{\ln (\delta)}{c_{2}\sqrt{\ln (\delta n)n}}\\
&\leq3(1+\epsilon)\E\paren{\sup_{f\in\F}\absj{(P_{n}-P)f}}+\paren{3\sqrt{2}+\frac1{c_{2}}+\frac3{c_{2}\epsilon}}\sqrt{\frac{\ln(\delta)}n}.
\end{align*}
We apply Lemma~\ref{lemma:EsperanceSupVdisjoints} to the sets $A_x=x(V)$ and the real numbers $\alpha_x=\paren{P(x(V/\{i\}))\vee c_{2}^{2}\ln (\delta n) n^{-1}}^{-1/2}$. We have
$$\alpha^*\leq \sqrt{n(c_{2}^{2}\ln (\delta n))^{-1}}\,\;p^*=\sup_{x\in\x(G)}(P(x(V/\{i\})))^{-1}P(x(V))\leq 1.$$
Hence, 
\begin{align*}
\E&\paren{\sup_{f\in\F}\absj{(P_{n}-P)f}}\leq\frac{64}{\sqrt{n}}\sqrt{\ln\left(4\sqrt{n(c_{2}^{2}\ln (\delta n))^{-1}}\right)}\\
&+2048\frac{\ln\left(4\sqrt{n(c_{2}^{2}\ln (\delta n))^{-1}}\right)}{c_{2}\sqrt{n(\ln (\delta n))^{-1}}}\leq \paren{32\sqrt{2}+\frac{2048}{c_2}}\sqrt{\frac{\ln(n)}{n}}.
\end{align*}
Thus, for all $\epsilon>0$, with probability larger than $1-\delta^{-1}$, we have
\begin{align*}
\sup_{x\in\x'_{n}}&\sqrt{\Ph(x(V/\{i\}))}|\Ph_{i|V}(x)-P_{i|V}(x)|\\
&\leq 2\paren{\paren{99\sqrt{2}+\frac{6144}{c_2}}(1+\epsilon)+\frac1{c_2}\paren{1+\frac3{\epsilon}}}\sqrt{\frac{\ln(\delta n)}{n}}.
\end{align*}
We take $\epsilon=0.001$ to conclude the proof.

\subsection{Proof of Theorem \ref{theo:GenModSel}:}
It comes from Theorem \ref{theo:controlVar} that, for all subsets $V$ in $\G_n$, we have,
$$P\paren{\left\|\Ph_{i|V}-P_{i|V}\right\|_{\infty}\leq c_2\sqrt{\frac{\ln(N_{n}\delta n)}{n\ph_-^V}}}\geq 1-\frac1{N_{n}\delta}.$$
We use a union bound to get that,
$$P\paren{\forall V\in\G_{n},\;\left\|\Ph_{i|V}-P_{i|V}\right\|_{\infty}\leq c_2\sqrt{\frac{\ln(N_{n}\delta n)}{n\ph_-^V}}}\geq 1-\delta^{-1}.$$
Hereafter in the proof of Theorem \ref{theo:GenModSel}, we denote by $v_n^{V}=\sqrt{\ln(N_{n}\delta n)\paren{n\ph_-^V}^{-1}}$ and by
$$\Omega=\left\{\forall V\in\G_n,\; \left\|\Ph_{i|V}-P_{i|V}\right\|_{\infty}\leq c_2v^{V}_n\right\}.$$
We have proved that $P(\Omega)\geq 1-\delta^{-1}$. Let $C>c_2$ and denote, for short $\Gh=\Gh(C,\delta,\G_{n})$. By definition of $\Gh$, for all $V\in \G_{n}$,
$$\left\|P_{i|G}\right\|_{\infty}-\left\|\Ph_{i|\Gh}\right\|_{\infty}+Cv^{\Gh}_n\leq \left\|P_{i|G}\right\|_{\infty}-\left\|\Ph_{i|V}\right\|_{\infty}+C\pen(V).$$
Hence, on $\Omega$, for all $V$ in $\G_n$,
\begin{align}\label{eq:interm1}
\left\|P_{i|G}\right\|_{\infty}-\left\|P_{i|\Gh}\right\|_{\infty}+(C-c_2)v^{\Gh}_n\leq \left\|P_{i|G}\right\|_{\infty}-\left\|P_{i|V}\right\|_{\infty}+(C+c_2)\pen(V).
\end{align}
From Assumption {\bf H1}, $\left\|P_{i|G}\right\|_{\infty}-\left\|P_{i|\Gh}\right\|_{\infty}\geq \kappa_{\min}\left\|P_{i|G}-P_{i|\Gh}\right\|_{\infty}$ and from the triangular inequality, $, \;\left\|P_{i|G}\right\|_{\infty}-\left\|P_{i|V}\right\|_{\infty}\leq \left\|P_{i|G}-P_{i|V}\right\|_{\infty}.$
Plugging these inequalities in (\ref{eq:interm1}), we obtain that, for all $V\in \G_{n}$,
\begin{align}\label{eq:interm2}
\kappa_{\min}\left\|P_{i|G}-P_{i|\Gh}\right\|_{\infty}+(C-c_2)v^{\Gh}_n\leq \left\|P_{i|G}-P_{i|V}\right\|_{\infty}+(C+c_2)\pen(V).
\end{align}
On $\Omega$, for all $V\in \G_{n}$, we have then
\begin{align*}
\left\|\Ph_{i|\Gh}-P_{i|G}\right\|_{\infty}&\leq\left\|\Ph_{i|\Gh}-P_{i|\Gh}\right\|_{\infty}+ \left\|P_{i|\Gh}-P_{i|G}\right\|_{\infty}\leq c_2v^{\Gh}_n+ \left\|P_{i|G}-P_{i|\Gh}\right\|_{\infty}\\
&\leq \max\left(\frac1{\kappa_{\min}},\frac{c_2}{C-c_2}\right)\left(\kappa_{\min}\left\|P_{i|G}-P_{i|\Gh}\right\|_{\infty}+(C-c_2)v^{\Gh}_n\right)
\end{align*}

\begin{align*}
\left\|\Ph_{i|\Gh}-P_{i|G}\right\|_{\infty}&\leq \max\left(\frac1{\kappa_{\min}},\frac{c_2} {C-c_2}\right)\left(\left\|P_{i|G}-P_{i|V}\right\|_{\infty}+(C+c_2)\pen(V)\right)\\
&\leq K(c_2,C,\kappa_{\min})\left(\left\|P_{i|G}-P_{i|V}\right\|_{\infty}+\pen(V)\right).
\end{align*}
\subsection{Proof of Corollary \ref{coro:ONProblem}:}
It comes from \cite{Ma07} Proposition 2.5 p 20 that
$$N_{m,M}=|\G_{m,M}|=\sum_{k=0}^mC_M^k\leq \left(\frac{eM}{m}\right)^m\leq M^m\;{\rm hence}\; \ln(N_{m,M})\leq m\ln(M).$$
Hence, from Theorem \ref{theo:GenModSel}, with probability larger than $1-\delta^{-1}$, we have
\begin{equation}\label{eq:ModSelGibbs}
\left\|\Ph_{i|\Gh_{\delta}(C))}-P_{i|G}\right\|_{\infty}\leq K\inf_{V\in \G_{m,M}}\left\{\left\|P_{i|G}-P_{i|V}\right\|_{\infty}+\sqrt{\frac{\ln(nM^{m}\delta)}{n\ph^V_-}}\right\}.
\end{equation}
For all $|V|> \log_{2}(n)$, there is at least one configuration in $\x(V)$ that is not observed, hence $\ph^V_-=1/n$. Therefore, for all $m\geq \log_{2}(n)$,
\begin{align*}
\inf_{V\in \G_{\log_{2}(n),M}}&\left\{\left\|P_{i|G}-P_{i|V}\right\|_{\infty}+\sqrt{\frac{\Gamma_{M}(\delta)}{n\ph^V_-}}\right\}= \inf_{V\in \G_{m,M}}\left\{\left\|P_{i|G}-P_{i|V}\right\|_{\infty}+\sqrt{\frac{\Gamma_{M}(\delta)}{n\ph^V_-}}\right\}.
\end{align*}
Taking $m=M$, (\ref{eq:ModSelGibbs}) yields the corollary.
\subsection{Proof of Theorem \ref{theo:IntGraphEst}:}
Let $\Omega$ be the event defined in the proof of Theorem~\ref{theo:GenModSel} for the collection $\G_{n}$ and let $\Gh=\Gh_{\delta}(C)$. We have $P(\Omega^{c})\leq \delta^{-1}$ and, on $\Omega$, from Corollary \ref{coro:ONProblem},
\begin{equation*}
\left\|P_{i|\Gh}-P_{i|G}\right\|_{\infty}\leq K\inf_{V\subset V_M}\set{\left\|P_{i|G}-P_{i|V}\right\|_{\infty}+v_{n}^{V}(\delta)}.
\end{equation*}
By definition of $\Psi$, denoting by $l_n=\sqrt{n^{-1}\Gamma_{M}(\delta)}$, we have
$$\inf_{V\subset V_M}\set{\left\|P_{i|G}-P_{i|V}\right\|_{\infty}+v_{n}^{V}(\delta)}=\inf_{v>0}\set{\Psi(v)+vl_n}.$$
Let $v^*$ be the smallest solution of the equation $vl_n=\Psi(v)$. As $\Psi$ is non-increasing and $v\mapsto vl_n$ is non decreasing, we have
$$ \Psi(v^*)\leq \inf_{v>0}\set{\Psi(v)+vl_n}\leq 2\Psi(v^*).$$
Thus, on $\Omega$, we have
\begin{equation*}
\left\|P_{i|\Gh}-P_{i|G}\right\|_{\infty}\leq 2K\Psi(v^*).
\end{equation*}
Let $\Omega_2$ be the event defined in {\bf H2}. Let $r<1$ and $\omega$ in $\Omega^*=\Omega \cap \Omega_{2}$ such that $\ph_-^{\Gh}(\omega)\geq (rv^*)^{-2}.$ From assumption {\bf H2} applied to $v=rv^*$, $K=r^{-1}$,
$$2K\Psi(v^*)\geq \left\|P_{i|\Gh}(\omega)-P_{i|G}\right\|_{\infty}\geq \Psi(rv^*)\geq C_{\Psi}^{-1}r^{-\alpha}\Psi(v^*).$$
Hence $r\geq (2C_{\Psi}K)^{-1/\alpha}$. Thus, on $\Omega^*$, we have
$$\left\|P_{i|\Gh}-P_{i|G}\right\|_{\infty}\leq 2Kl_nv^*\leq C_{\Psi}^{-1/\alpha}(2K)^{1-1/\alpha}v_n^{\Gh}.$$
By the triangular inequality, we have
\begin{align*}
\sup_{x\in\x(G)}|(\Ph_{i|V}(x)-\Ph_{i|V}(x_j))&-(P_{i|G}(x)-P_{i|G}(x_j))|\\
&\leq 2\left(\norm{\Ph_{i|V}-P_{i|V}}_{\infty}+\left\|P_{i|V}-P_{i|G}\right\|_{\infty}\right).
\end{align*}
Hence, on $\Omega^*$,
\begin{align*}
\absj{\omega^{\Gh}_{i,j}(\Ph)- \omega^{G}_{i,j}(P)}&\leq 2\left(c_{2}v^{\Gh}_n+\left\|P_{i|\Gh}-P_{i|G}\right\|_{\infty}\right)\\
&\leq 2\paren{c_{2}+C_{\Psi}^{-1/\alpha_{\Psi}}(2K)^{1-1/\alpha_{\Psi}}}v_n^{\Gh}.
\end{align*}
Let $c_{\infty}=2\paren{c_{2}+C_{\Psi}^{-1/\alpha_{\Psi}}(2K)^{1-1/\alpha_{\Psi}}}$. It comes from this last inequality that, on $\Omega^{*}$,
\begin{align*}
&\set{j\in V_{M},\;\omega^{G}_{i,j}(P)\geq (c+c_{\infty})v_{n}^{\Gh}}\\
&\subset \Gh_{i}^{\Gh}(c)\subset \set{j\in V_{M},\;\omega^{G}_{i,j}(P)\geq (c-c_{\infty})v_{n}^{\Gh}}.
\end{align*}

\subsection{Proof of Theorem \ref{theo:controlBiasGibbs}:}
In all the proof, for all subsets $V$, $V'$ of $G$ such that $V\cap V'=\emptyset$, for all $(x,y)$ in $\x(V)\times \x(V')$, let $x(V)\oplus y(V')$ be the configuration on $\x(V\cup V')$ such that, for all $j$ in $V$ $x(V)\oplus y(V')(j)=x(j)$ and for all $j$ in $V'$, $x(V)\oplus y(V')(j)=y(j)$. Let $V$ be a finite subset of $G$ and let $x$ be a configuration on $\x(G)$.
\begin{align}\label{eq:gibbsH2}
P_{i|G}(x)-&P_{i|V}(x)=\int (P_{i|G}(x)-P_{i|G}(x(V)\oplus y(G/V)))dP(y(G/V)|x(V/\{i\}))
\end{align}
From the definition of a Gibbs measure, we have
\begin{align}\label{eq:PiG=fij}
&P_{i|G}(x)-P_{i|G}(x(V)\oplus y(G/V))\nonumber\\
&=\frac{e^{-\sum_{j\in G}g_{i,j}(x(i),x(j))}\left(e^{\sum_{j\notin V}(g_{i,j}(x(i),x(j))-g_{i,j}(x(i),y(j)))}-1\right)}{\left(1+e^{-\sum_{j\in G}g_{i,j}(x(i),x(V)\oplus y(G/V)(j))}\right)\left(1+e^{-\sum_{j\in G}g_{i,j}(x(i),x(j))}\right)}
\end{align}
Hence,
\begin{align*}
|P_{i|G}(x)&-P_{i|G}(x(V)\oplus y(G/V))|\\
&\leq \frac{e^{2r}}{(1+e^{-2r})^2}\left|e^{\sum_{j\notin V}(g_{i,j}(x(i),x(j))-g_{i,j}(x(i),y(j)))}-1\right|.
\end{align*}
Let us now give the following lemma, whose proof is immediate from the convexity of $x\mapsto e^x$.
\begin{lemma}\label{lem:Ana2}
For all real numbers $r>0$, for all $x$ in $[-4r,4r]$, we have
$$\frac{1-e^{-4r}}{4r}|x|\leq |e^x-1|\leq\frac{e^{4r}-1}{4r}|x|.$$
\end{lemma}

\noindent
We deduce from Lemma \ref{lem:Ana2} that
\begin{align*}
|P_{i|G}(x)&-P_{i|G}(x(V)\oplus y(G/V))|\\
&\leq \frac{(e^{4r}-1)e^{2r}}{4r(1+e^{-2r})^2}\left|\sum_{j\notin V}(g_{i,j}(x(i),x(j))-g_{i,j}(x(i),y(j)))\right|.\nonumber
\end{align*}
It is clear that, for all $x$ in $\x(G)$,
\begin{align*}
\int\left|\sum_{j\notin V}(g_{i,j}(x(i),x(j))-g_{i,j}(x(i),y(j)))\right|&dP(y(G/V)|x(V/\{i\}))\\
&\leq \sum_{j\notin V}\omega_{i,j}(f)P(y(j)\neq x(j)|x(V/\{i\})).
\end{align*}
The upper bound comes then from the inequality $P(y(j)\neq x(j)|x(V/\{i\}))\leq 1$.

For the lower bound, let, for all $a$ in $A$, $x^a_{\max}$ be the configuration such that, for all $j$ in $G$, $g_{i,j}(a,x_{\max}^a(j))=\left\|g_{i,j}^a\right\|$ and let $x^a_{\min}$ be the configuration such that, for all $j$ in $G$, $g_{i,j}(a,x^a_{\min}(j))=\inf_{b\in A}g_{i,j}(a,b)$. From Lemma \ref{lem:Ana2}, we have
\begin{align}\label{eq:intfij2}
e^{\sum_{j\notin V}\left\|g_{i,j}^{x(i)}\right\|-g_{i,j}(x(i),y(j))}-1&=\left|e^{\sum_{j\notin V}\left\|g_{i,j}^{x(i)}\right\|-g_{i,j}(x(i),y(j))}-1\right|\nonumber\\
&\geq \frac{1-e^{-4r}}{4r}\sum_{j\notin V}\left\|g_{i,j}^{x(i)}\right\|-g_{i,j}(x(i),y(j)).
\end{align}
Finally, we have
\begin{align}\label{eq:intfij}
\int \left\|g_{i,j}^{x(i)}\right\|&-g_{i,j}(x(i),y(j))dP(y(G/V)|x(V/\{i\}))\nonumber\\
&\geq P(y(j)=x_{\min}^{x(i)}(j)|x(V/\{i\}))\omega_{i,j}(f)\geq \frac{\omega_{i,j}(f)}{1+e^{2r}}.
\end{align}
Using successively inequalities (\ref{eq:gibbsH2}), (\ref{eq:PiG=fij}), (\ref{eq:intfij2}) and (\ref{eq:intfij}) with $x=x^{x(i)}_{\max}$, we obtain

\begin{align*}
&\sup_{x\in \x(G)}P_{i|G}(x)-P_{i|V}(x)\\
&\geq \int (P_{i|G}(x^{x(i)}_{\max})-P_{i|G}(x^{x(i)}_{\max}(V)\oplus y(G/V)))dP(y(G/V)|x(V/\{i\}))\\
&=\frac{e^{-2r}}{1+e^{-2r}}\int\frac{e^{\sum_{j\notin V}\left\|g_{i,j}^{x(i)}\right\|-g_{i,j}(x(i),y(j))}-1}{1+e^{-\sum_{j\in V}\left\|g_{i,j}^{x(i)}\right\|-\sum_{j\notin V}g_{i,j}(x(i),y(j))}}dP(y(G/V)|x(V/\{i\}))\\
&\geq \frac{e^{-2r}}{(1+e^{2r})^2}\int \left(e^{\sum_{j\notin V}\left\|g_{i,j}^{x(i)}\right\|-g_{i,j}(x(i),y(j))}-1\right)dP(y(G/V)|x(V/\{i\}))
\end{align*}
Hence,
\begin{align*}
&\sup_{x\in \x(G)}P_{i|G}(x)-P_{i|V}(x)\\
&\geq \frac{(1-e^{-4r})e^{-2r}}{4r(1+e^{2r})^2}\int \left(\sum_{j\notin V}\left\|g_{i,j}^{x(i)}\right\|-g_{i,j}(x(i),y(j))\right)dP(y(G/V)|x(V/\{i\}))\\
&\geq \frac{(1-e^{-4r})e^{-2r}}{4r(1+e^{2r})^3}\sum_{j\notin V}\omega_{i,j}(f).
\end{align*}
Let us now check that $P$ satisfies assumption {\bf H1}. Let $x$ in $\x(V)$. Using successively inequalities (\ref{eq:gibbsH2}), (\ref{eq:PiG=fij}), (\ref{eq:intfij2}) and (\ref{eq:intfij}) with $x=x(V)\oplus x^{x(i)}_{\max}(G/V)$, we obtain, as in the previous proof,
\begin{align*}
P_{i|G}(x(V)\oplus x^{x(i)}_{\max}(G/V))-P_{i|V}(x)&\geq \frac{(1-e^{-4r})e^{-2r}}{4r(1+e^{2r})^3}\sum_{j\notin V}\omega_{i,j}(f).
\end{align*}
Taking $x$ such that $P_{i|V}(x)=\left\|P_{i|V}\right\|_{\infty}$ and using that $P_{i|G}(x(V)\oplus x^{x(i)}_{\max}(G/V))\leq \left\|P_{i|G}\right\|_{\infty}$, we obtain
$$\left\|P_{i|G}\right\|_{\infty}-\left\|P_{i|V}\right\|_{\infty}\geq \frac{(1-e^{-4r})e^{-2r}}{4r(1+e^{2r})^3}\sum_{j\notin V}\omega_{i,j}(f).$$
This yields the theorem thanks to inequality (\ref{eq:ControlBias1}).

\section{Proof of Theorem \ref{theo:BornesRisk}:}
\begin{align*}
&\absj{\ph(i,j)-\ph(i)\ph(j)-\paren{P(x_{1}(i,j))-P(x_{1}(i))P(x_{1}(j))}}\\
&\leq\absj{\Ph(x_{1}(i,j))-P(x_{1}(i,j))}+\absj{\Ph(x_{1}(j))-P(x_{1}(j))}+\absj{\Ph(x_{1}(i))-P(x_{1}(i))}
\end{align*}
We use Hoeffding's inequality (see for example \cite{Ma07} Proposition 2.7) to the functions $t=1_{x_{1}(i,j)},\;1_{x_{1}(i)},\;1_{x_{1}(j)}$, for all $x>0$, we have
$$P\paren{\absj{(P_{n}-P)t}>\sqrt{\frac x{2n}}}\leq 2e^{-x}.$$
Hence, a union bound gives that, on a set $\Omega(\delta)$ satisfying $P(\Omega(\delta)^{c})\leq \delta^{-1}$, for all $j$ in $V_{M}$
\begin{align*}
\absj{\Ph(x_{1}(i,j))-P(x_{1}(i,j))}&+\absj{\Ph(x_{1}(j))-P(x_{1}(j))}\\
&+\absj{\Ph(x_{1}(i))-P(x_{1}(i))}\leq 3\sqrt{\frac{\ln(6M\delta)}{2n}}.
\end{align*}
Moreover, we have
\begin{align}
&P_{i|j}(x_{1})-P(x_{1}(i))=\paren{P_{i|j}(x_{1})-P_{i|\emptyset}(x_{1})}\label{eq:Cov1}=\\
&\int P_{i|G}(x_{1}(i,j)\oplus y(G/\{i,j\}))-P_{i|G}(x_{1}(i)\oplus y(G/\{i\}))dP(y(G/\{i\})|x_{1}(i)).\nonumber
\end{align}
From the definition of a Gibbs measure, we have
\begin{align}
P_{i|G}&(x_{1}(i,j)\oplus y(G/\{i,j\}))-P_{i|G}(x_{1}(i)\oplus y(G/\{i\}))\label{eq:Cov2}\\
&=\frac{e^{-\sum_{j\in G}g_{i,j}(x_{1}(i),y(j))}\left(e^{(g_{i,j}(x_{1}(i),x_{1}(j))-g_{i,j}(x_{1}(i),y(j)))}-1\right)}{\left(1+e^{-\sum_{j\in G}g_{i,j}(x_{1}(i),x_{1}(i,j)\oplus y(G/\{i,j\}))}\right)\left(1+e^{-\sum_{j\in G}g_{i,j}(x_{1}(i),y(j))}\right)}\nonumber.
\end{align}
We can assume that, without loss of generality that 
$g_{i,j}(x_{1}(i),x_{1}(j))=\norm{g_{i,j}}$ and therefore that 
$$e^{g_{i,j}(x_{1}(i),x_{1}(j))-g_{i,j}(x_{1}(i),y(j))}-1=\absj{e^{g_{i,j}(x_{1}(i),x_{1}(j))-g_{i,j}(x_{1}(i),y(j))}-1}.$$
It comes then from Lemma \ref{lem:Ana2} that
\begin{align*}
\frac{1-e^{-4r}}{4r}&\absj{g_{i,j}(x_{1}(i),x_{1}(j))-g_{i,j}(x_{1}(i),y(j))}\\
&\leq e^{g_{i,j}(x_{1}(i),x_{1}(j))-g_{i,j}(x_{1}(i),y(j))}-1\\
&\leq \frac{e^{4r}-1}{4r}\absj{g_{i,j}(x_{1}(i),x_{1}(j))-g_{i,j}(x_{1}(i),y(j))}
\end{align*}
Hence
\begin{align}\label{eq:Cov3}
\frac{1-e^{-4r}}{4r}&\absj{\omega_{i,j}(f)}1_{y(j)\neq x_{1}(j)}\nonumber\\
&\leq e^{g_{i,j}(x_{1}(i),x_{1}(j))-g_{i,j}(x_{1}(i),y(j))}-1\leq \frac{e^{4r}-1}{4r}\absj{\omega_{i,j}(f)}
\end{align}
Using successively (\ref{eq:Cov1}), (\ref{eq:Cov2}), (\ref{eq:Cov3}), we deduce that
\begin{align*}
\frac{e^{-2r}(1-e^{-4r})}{4r(1+e^{2r})^{3}}&\absj{\omega_{i,j}(f)}\leq \frac{e^{-2r}(1-e^{-4r})}{4r(1+e^{2r})^{2}}P(y(j)\neq x_{1}(j)|x_{1}(i))\absj{\omega_{i,j}(f)}\\
&\leq \absj{P(x_{1}(i,j))-P(x_{1}(i))P(x_{1}(j))} \leq \frac{e^{2r}(e^{4r}-1)}{4r(1+e^{-2r})^{2}}\absj{\omega_{i,j}(f)}.
\end{align*}
We conclude that, on $\Omega(\delta)$,
\begin{align*}
&\frac{e^{-2r}(1-e^{-4r})}{4r(1+e^{2r})^{3}}\absj{\omega_{i,j}(f)}-3\sqrt{\frac{\ln(6M\delta)}{2n}}\\
&\leq \absj{\Ph(x_{1}(i,j))-\Ph(x_{1}(i))\Ph(x_{1}(j))}\leq \frac{e^{2r}(e^{4r}-1)}{4r(1+e^{-2r})^{2}}\absj{\omega_{i,j}(f)}+3\sqrt{\frac{\ln(6M\delta)}{2n}}
\end{align*}
All the sites $j$ such that 
\begin{equation}\label{eq:thoseinside}
\absj{\omega_{i,j}(f)}\geq \frac{4r(1+e^{2r})^{3}}{e^{-2r}(1-e^{-4r})}\paren{\eta+3\sqrt{\frac{\ln(6M\delta)}{2n}}}
\end{equation}
belong to $\Vh(\eta)$. All the sites such that
\begin{equation}\label{eq:thoseoutside}
\absj{\omega_{i,j}(f)}< \frac{4r(1+e^{-2r})^{2}}{e^{2r}(e^{4r}-1)}\paren{\eta-3\sqrt{\frac{\ln(6M\delta)}{2n}}}
\end{equation}
do not belong to $\Vh(\eta)$.\\
We use then Theorem \ref{theo:GenModSel} with the collection $\G=\set{V\subset V(\eta_{ms})}$. Its cardinality is bounded by $n^{\kappa}$.  There exist a constant $K$ and an event $\Omega_{2}(\delta)$, with probability $1-\delta$, such that, on $\Omega_{2}(\delta)$,
$$\norm{\Ph_{i|\Gh}-P_{i|G}}\leq K\inf_{V\in\G}\set{\norm{P_{i|V}-P_{i|G}}+\sqrt{\frac{\ln(n^{\kappa}\delta)}{n\ph_{-}^{V}}}}.$$
We use Theorem \ref{theo:controlBiasGibbs} to say that 
$$\norm{P_{i|V}-P_{i|G}}\leq C_{r}\sum_{j\notin V}\absj{\omega_{i,j}(f)}=C_{r}\paren{\sum_{j\notin \Vh(\eta)}\absj{\omega_{i,j}(f)}+\sum_{j\in \Vh(\eta)/V}\absj{\omega_{i,j}(f)}}.$$
We deduce from (\ref{eq:thoseinside}) that, on $\Omega(\delta)$,
$$\sum_{j\notin \Vh(\eta)}\absj{\omega_{i,j}(f)}\leq \sum_{j\in V(\eta,\delta,M)}\absj{\omega_{i,j}(f)}.$$
We choose $\Omega_{*}(\delta)=\Omega(\delta)\cap \Omega_{2}(\delta)$ to conclude the proof.

\section{Appendix}\label{section:Appendix}
In this Appendix, we recall the bound given by Bousquet \cite{Bo02} for the deviation of the supremum of the empirical process.
\begin{theorem}
Let $X_1,...,X_n$ be i.i.d. random variables valued in a measurable space $(A,\x)$. Let $\F$ be a class of real valued functions, defined on $A$ and bounded by $b$. Let $v^2=\sup_{f\in \F}P[(f-Pf)^2]$ and $Z=\sup_{f\in \F}(P_n-P)f$. Then, for all $x>0$,
\begin{equation}
P\left(Z>\E(Z)+\sqrt{\frac2n(v^2+2b\E(Z))x}+\frac{bx}{3n}\right)\leq e^{-x}.
\end{equation}
\end{theorem}

\noindent
Let us recall some well known tools of empirical processes theory.
\begin{definition}
The covering number $N(\epsilon,T,d)$ is the minimal number of balls of radius $\epsilon$ with centers in $T$ needed to cover $T$. The entropy is the logarithm of the covering number $H(\epsilon,T,d)=\ln (N(\epsilon,T,d)).$
\end{definition}
\begin{definition}
An $\epsilon$-separated subset of $T$ is a subset $\{t_k\}$ of elements of $T$ whose pairwise distance is strictly larger than $\epsilon$. The packing number $M(\epsilon,T,d)$ is the maximum size of an $\epsilon$-separated subset of $T$.
\end{definition}

\noindent
Those quantities are related by the famous following lemma.
\begin{lemma}
(Kolmogorov and Tikhomirov \cite{KT63}) Let $(T,d)$ be a metric space and let $\epsilon>0$,
$$N(\epsilon,T,d)\leq M(\epsilon,T,d)\leq N(\epsilon/2,T,d).$$
\end{lemma}

The following result can be derived from classical chaining arguments (see for example \cite{Bo02}).
\begin{lemma}\label{lem:Esup}
Let $\F$ be a class of functions, let $d_{2,P_n}(t,t')=\sqrt{P_n[(t-t')^2]}$ and $D_n=\sqrt{\sup_{t\in \F}P_n(t^2)}$ then 
$$\E\left(\sup_{t\in \F}|(P_n-P)t|\right)\leq \frac{16\sqrt{2}}{\sqrt{n}}\E\left(\int_{0}^{D_n/2}H^{1/2}(u,\F,d_{2,P_n})du\right).$$
\end{lemma}

The next result was used to obtain our concentration inequalities.
\begin{lemma}\label{lemma:EsperanceSupVdisjoints}
Let $(A_i)_{i\in I}$ be a collection of sets such that, for all $i,j\in I$, $A_i\cap A_j=\emptyset$ and let $(\alpha_i)_{i\in I}$ be a collection of positive real numbers. Let $Z_I=\sup_{t\in \F_I}|(P_n-P)t|$, where $\F_I=\{t_i=\alpha_i1_{A_i}\}$ and $P_n$ is the empirical measure. Let $\alpha^*=\sup_{i\in I}\alpha_i,$ $p_*=\sup_{i\in I}\alpha_i^2P(A_i)$. We have
\begin{equation}\label{eq:Esup}
\E\left(\sup_{t\in \F_I}|(P_n-P)t|\right)\leq\frac{64}{\sqrt{n}}\sqrt{p^*\ln\left(\frac{4\alpha^*}{\sqrt{p^*}}\right)}+\frac{2048}n\alpha^*\ln\left(\frac{4\alpha^*}{\sqrt{p^*}}\right).
\end{equation}
\end{lemma}

\noindent 
In order to apply Lemma~\ref{lem:Esup} to $\F=\F_I$, we compute the entropy of $\F_I$. For all $i\neq j$, since $A_i\cap A_j=\emptyset$,
$$(t_i-t_j)^2=\left(\alpha_i1_{A_i}-\alpha_j1_{A_j}\right)^2=\alpha_i^21_{A_i}+\alpha_j^21_{A_j}.$$
Hence $d_{2,P_n}(t_i,t_j)=\sqrt{\alpha_i^2P_n(A_i)+\alpha_j^2P_n(A_j)}.$\\
Consider an $\epsilon$-separated set $T_{\epsilon}=\{t_{i_1},...,t_{i_N}\}$ in $(\F_I,d_{2,P_n})$ (see also the definition in the appendix), it comes from the previous computation that, for all $k\neq k'$,
$$\alpha_{i_k}^2P_n(A_{i_k})+\alpha_{i_{k'}}^2P_n(A_{i_{k'}})\geq \epsilon^2.$$
Hence, there is at least $N-1$ indexes $k\in\{1,...,N\}$ such that $\alpha_{i_k}^2P_n(A_{i_k})\geq \epsilon^2/2.$
It follows that
$$1=\sum_{i\in I}P_n(A_i)\geq \sum_{k=1}^NP_n(A_{i_k})\geq \frac{\epsilon^2(N-1)}{2(\alpha^*)^2}.$$
Hence $N\leq 1+2(\alpha^*)^2\epsilon^{-2}$, thus $H(\epsilon,\F_I,d_{2,P_n})\leq \ln\left(1+2(\alpha^*)^2\epsilon^{-2}\right).$
We deduce from this inequality and Lemma \ref{lem:Esup} that

\begin{align}
\E\left(\sup_{t\in \F_I}|(P_n-P)t|\right)&\leq \frac{16\sqrt{2}}{\sqrt{n}}\E\left(\int_0^{\sqrt{\hat{p}_n^*}/2}\sqrt{\ln\left(1+2(\alpha^*)^2\epsilon^{-2}\right)}d\epsilon\right)\notag\\
&\leq\frac{32}{\sqrt{n}}\E\left(\int_0^{\sqrt{\hat{p}_n^*}/2}\sqrt{\ln\left(2\alpha^*\epsilon^{-1}\right)}d\epsilon\right) ,\label{eq:expecsup}
\end{align}
where $\hat{p}_n^*=\sup_{i\in I}\alpha_i^2P_n(A_i)$. Now, let us recall the following elementary lemma.
\begin{lemma}\label{lem:intsqrtln}
For all positive real numbers $K,A$ such that $K/A>e$, we have
$$\int_0^{A}\sqrt{\ln (Kx^{-1})}dx\leq 2A\sqrt{\ln\left(\frac KA\right)}$$
\end{lemma}

\noindent
Actually,
\begin{align*}
\int_0^{A}\sqrt{\ln (Kx^{-1})}dx=K\int_{K/A}^{\infty}\frac{\sqrt{\ln (x)}}{x^2}dx=A\sqrt{\ln\left(\frac KA\right)}+\frac K2\int_{K/A}^{\infty}\frac{1}{u^2\sqrt{\ln u}}du.
\end{align*}
Since $K/A>e$, $\frac{1}{u^2\sqrt{\ln u}}\leq \frac{\sqrt{\ln u}}{u^2}$ on $[K/A,\infty[$. The result follows.\\
By definition, $\hat{p}_n\leq (\alpha^*)^2$, hence $2\alpha^*/(\sqrt{\hat{p}_n^*}/2)\geq 4>e$, we deduce from Lemma \ref{lem:intsqrtln} that
$$\E\left(\sup_{t\in \F_I}|(P_n-P)t|\right)\leq \frac{32}{\sqrt{n}}\E\left[\sqrt{\hat{p}_n^*}\sqrt{\ln\left(\frac{4\alpha^*}{\sqrt{\hat{p}_n^*}}\right)}\right].$$
Let us now give another simple lemma.
\begin{lemma}\label{lem:Analyse}
The function $f:x\mapsto x\sqrt{\ln(K/x)}$, defined on $(0,K)$ is positive, non decreasing on $(0,K/e^{1/2})$ and strictly concave.
\end{lemma}

\noindent
The proof of the lemma is straightforward from the computations
$$f'(x)=\sqrt{\ln(K/x)}-\frac1{2\sqrt{\ln(K/x)}},\;f"(x)=-\frac1{2x\sqrt{\ln(K/x)}}-\frac1{4x(\sqrt{\ln(K/x)})^3}.$$
Applying Lemmas \ref{lem:Analyse}, \ref{lem:intsqrtln}, and Jensen's inequality to the right hand side of (\ref{eq:expecsup}) we have that
$$\E\left(\sup_{t\in \F_I}|(P_n-P)t|\right)\leq \frac{32}{\sqrt{n}}\E\left(\sqrt{\hat{p}_n^*}\right)\sqrt{\ln\left(\frac{4\alpha^*}{\E\left(\sqrt{\hat{p}_n^*}\right)}\right)}.$$
Now it comes from Jensen inequality that 
$$\E\left[\sqrt{\hat{p}_n^*}\right]\leq \sqrt{\E\left[\hat{p}_n^*\right]}\leq \sqrt{p^*}+\sqrt{\alpha^*\E\left(\sup_{t\in \F_I}|(P_n-P)t|\right)}.$$
It is clear from its definition that $p^*\leq (\alpha^*)^2$. Moreover, as $P_n$ and $P$ are probability measures, we have, for all $t$ in $\F_I$, $|(P_n-P)t|\leq 2\alpha^*$. Hence, $\sqrt{\alpha^*\E\left(\sup_{t\in \F_I}|(P_n-P)t|\right)}\leq \sqrt{2}\alpha^*$. We deduce from these inequalities that
$$\sqrt{p^*}+\sqrt{\alpha^*\E\left(\sup_{t\in \F_I}|(P_n-P)t|\right)}\leq (1+\sqrt{2})\alpha^*\leq (4\alpha^*)/e^{1/2}.$$
Hence, it comes from Lemma \ref{lem:Analyse} that, if $E=\E\left(\sup_{t\in \F_I}|(P_n-P)t|\right)$
\begin{align*}
E&\leq \frac{32}{\sqrt{n}}\left(\sqrt{p^*}+\sqrt{\alpha^*E}\right)\sqrt{\ln\left(\frac{4\alpha^*}{\sqrt{p^*}+\sqrt{\alpha^*E}}\right)}\\
&\leq \frac{32}{\sqrt{n}}\left(\sqrt{p^*}+\sqrt{\alpha^*E}\right)\sqrt{\ln\left(\frac{4\alpha^*}{\sqrt{p^*}}\right)}.
\end{align*}
It is then straightforward that (\ref{eq:Esup}) holds.

\bibliographystyle{imsart-nameyear}
\bibliography{bibliolerasle}

\end{document}